\documentclass[11pt]{article}

\usepackage{fullpage,times,url}
\RequirePackage[authoryear]{natbib}

\usepackage{amsthm,amsfonts,amsmath,amssymb,epsfig,color,float,graphicx,verbatim,mathtools}
\usepackage{algorithm,algorithmic}

\usepackage{hyperref}
\hypersetup{
	colorlinks   = true, %Colours links instead of ugly boxes
	urlcolor     = blue, %Colour for external hyperlinks
	linkcolor    = blue, %Colour of internal links
	citecolor   = red %Colour of citations
}

\newtheorem{theorem}{Theorem}

\newtheorem{lemma}{Lemma}

\newtheorem{assumption}{Assumption}

\newcommand{\reals}{\mathbb{R}}
\newcommand{\E}{\mathbb{E}}

\newcommand{\bi}{\mathbf{i}}
\newcommand{\bj}{\mathbf{j}}
\newcommand{\be}{\mathbf{e}}
\newcommand{\bx}{\mathbf{x}}
\newcommand{\bw}{\mathbf{w}}

\newcommand{\bu}{\mathbf{u}}
\newcommand{\bv}{\mathbf{v}}

\newcommand{\Ucal}{\mathcal{U}}
\newcommand{\Ical}{\mathcal{I}}
\newcommand{\Ocal}{\mathcal{O}}

\newcommand{\Wcal}{\mathcal{W}}

\newcommand{\norm}[1]{\|#1\|}
\newcommand{\inner}[1]{\langle#1\rangle}

\newcommand{\secref}[1]{Sec.~\ref{#1}}

\renewcommand{\eqref}[1]{Eq.~(\ref{#1})}
\newcommand{\lemref}[1]{Lemma~\ref{#1}}

\newcommand{\thmref}[1]{Thm.~\ref{#1}}

\title{Oracle Complexity of Second-Order Methods\\ for Finite-Sum Problems}
\author{Yossi Arjevani and Ohad Shamir\\Department of Computer Science and Applied Mathematics\\
Weizmann Institute of Science}
\date{}

\begin{document}

\maketitle	
\begin{abstract}
	Finite-sum optimization problems are ubiquitous in machine learning, and 
	are commonly solved using first-order methods which rely on gradient 
	computations. Recently, there has been growing interest in 
	\emph{second-order} methods, which rely on both gradients and Hessians. In 
	principle, second-order methods can require much fewer iterations than 
	first-order methods, and hold the promise for more efficient algorithms. 
	Although computing and manipulating Hessians is prohibitive for 
	high-dimensional problems in general, the Hessians of individual functions 
	in finite-sum problems can often be efficiently computed, e.g. because they 
	possess a low-rank structure. Can second-order information indeed be used 
	to solve such problems more efficiently? In this paper, we provide evidence 
	that the answer -- perhaps surprisingly -- is negative, at least in terms 
	of worst-case guarantees. We also discuss what additional 
	assumptions and algorithmic approaches might potentially circumvent this 
	negative result. 
\end{abstract}

\section{Introduction}

We consider finite-sum problems of the form
\begin{equation}\label{eq:obj}
\min_{\bw\in\Wcal} F(\bw) = \frac{1}{n}\sum_{i=1}^{n}f_i(\bw),
\end{equation}
where $\Wcal$ is a closed convex subset of some Euclidean or Hilbert space, each $f_i$ is convex and $\mu$-smooth, and $F$ is $\lambda$-strongly convex\footnote{For a twice-differentiable function $f$, it is $\mu$-smooth and $\lambda$-strongly convex 
if $ \lambda I \preceq \nabla^2 f(\bw)\preceq \mu I$ for all $\bw\in\Wcal$.}. Such problems are ubiquitous in machine learning, for example in order to perform empirical risk minimization using convex losses.

To study the complexity of this and other optimization problems, it is common to consider an oracle model, where the optimization algorithm has no a-priori information about the objective function, and obtains information from an oracle which provides values and derivatives of the function at various domain points \citep{YudNem83}. The complexity of the algorithm is measured in terms of the number of oracle calls required to optimize the function to within some prescribed accuracy.

Existing lower bounds for finite-sum problems   show that using a first-order oracle, which given a point $\bw$ and index $i=1,\ldots,n$ returns $f_i(\bw)$ and $\nabla f_i(\bw)$, the number of oracle queries required to find an $\epsilon$-optimal solution is at least of order
\[
\Omega\left(n+\sqrt{\frac{n\mu}{\lambda}}\log\left(\frac{1}{\epsilon}\right)\right),
\]
either under algorithmic assumptions or assuming the dimension is sufficiently large\footnote{Depending on how $\epsilon$-optimality is defined precisely, and where the algorithm is assumed to start, these bounds may have additional factors inside the log. For simplicity, we present the existing bounds assuming $\epsilon$ is sufficiently small, so that a $\log(1/\epsilon)$ term dominates.} \citep{agarwal2014lower,lan2015optimal,woodworth2016tight,arjevani2016dimension}. This is matched (up to log factors) by existing approaches, and cannot be improved in general. 

An alternative to first-order methods are  \emph{second-order} methods, which also utilize Hessian information. A prototypical example is the Newton method, which given a (single) function $F$, performs iterations of the form
\begin{equation}\label{eq:newton}
\bw_{t+1} = \bw_{t}-\alpha_t \left(\nabla^2 F(\bw)\right)^{-1}\nabla F(\bw),
\end{equation}
where $\nabla F(\bw),\nabla^2 F(\bw)$ are the gradient and the Hessian of $F$ at $\bw$, and $\alpha_t$ is a step size parameter.  Second-order methods can have extremely fast convergence, better than those of first-order methods (i.e. quadratic instead of linear). Moreover, they can be invariant to affine transformations of the objective function, and provably independent of its strong convexity and smoothness parameters (assuming e.g. self-concordance) \citep{boyd2004convex}. A drawback of these methods, however, is that they can be computationally prohibitive. In the context of machine learning, we are often interested in high-dimensional problems (where the dimension $d$ is very large), and the Hessians are $d\times d$ matrices which in general may not even fit into computer memory. However, for optimization problems as in \eqref{eq:obj}, the Hessians of individual $f_i$ often have a special structure. For example, a very common special case of finite-sum problems in machine learning is empirical risk minimization for linear predictors, where $$f_i(\bw)=\ell_i(\inner{\bw,\bx_i}),$$ where $\bx_i$ is a training instance and $\ell_i$ is some loss function. In that case, assuming $\ell_i$ is twice-differentiable, the Hessian has the rank-1 form $\ell''_i(\inner{\bw,\bx_i})\bx_i\bx_i^\top$. Therefore, the memory and computational effort involved with storing and manipulating the Hessian of this function is merely linear (rather than quadratic) in $d$. Thus, it is tractable even for high-dimensional problems. 

Building on this, several recent papers proposed and analyzed second-order methods for finite-sum problems, which utilize Hessians of the individual functions $f_i$ (see for instance \citet{erdogdu2015convergence,agarwal2016second,pilanci2015newton,roosta2016sub,roosta2016sub2,bollapragada2016exact,xu2016sub} and references therein). These can all be viewed as approximate Newton methods, which replace the actual Hessian $\nabla^2 F(\bw)=\frac{1}{n}\sum_{i=1}^{n}\nabla^2 f_i(\bw)$ in \eqref{eq:newton} by some approximation, based for instance on the Hessians of a few individual functions $f_i$ sampled at random. One may hope that such methods can inherit the favorable properties of second-order methods, and improve on the performance of commonly used first-order methods.  

In this paper, we consider the opposite direction, and study \emph{lower} bounds on the number of iterations required by algorithms using second-order (or possibly even higher-order) information, focusing on finite-sum problems which are strongly-convex and smooth. We make the following contributions:
\begin{itemize}
\item First, as a more minor contribution, we prove that in the standard setting of optimizing a \emph{single} smooth and strongly convex function, second-order information cannot improve the oracle complexity compared to first-order methods (at least in high dimensions). Although this may seem unexpected at first, the reason is that the smoothness constraint must be extended to higher-order derivatives, in order for higher-order information to be useful. We note that this observation in itself is not new, and is briefly mentioned (without proof) in \citet[Section 7.2.6]{YudNem83}. Our contribution here is in providing a clean, explicit statement and proof of this result.
\item We then turn to present our main results, which state (perhaps surprisingly) that under some mild algorithmic assumptions, and if the dimension is sufficiently large, the oracle complexity of second-order methods for finite-sum problems is no better than first-order methods, \emph{even} if the finite-sum problem is composed of quadratics (which are trivially smooth to any order).
\item Despite this pessimistic conclusion, our results also indicate what 
assumptions and algorithmic approaches might be helpful in circumventing it. In 
particular, it appears that better, dimension-dependent performance may be 
possible, if the dimension is moderate and the $n$ individual functions in 
\eqref{eq:obj} are accessed adaptively, in a manner depending on the functions 
rather than fixed in advance (e.g. sampling them from a non-uniform 
distribution depending on their Hessians, as opposed to sampling them uniformly 
at random). This provides evidence to the necessity of adaptive sampling 
schemes, and a dimension-dependent analysis, which indeed accords with some 
recently proposed algorithms and derivations, e.g. 
\cite{agarwal2016second,xu2016sub}. We note that the limitations arising from 
oblivious optimization schemes (in a somewhat stronger sense) was also explored 
in \cite{arjevani2016dimension,arjevani2016iteration}. 
\end{itemize}

The paper is structured as follows: We begin in \secref{sec:single} with a 
lower bound for algorithms utilizing second-order information, in the simpler 
setting where there is a single function $F$ to be optimized, rather than a 
finite-sum problem. We then turn to provide our main lower bounds in 
\secref{sec:main}, and discuss their applicability to some existing approaches 
in \secref{sec:related}. We conclude in \secref{sec:summary}, where we also 
discuss possible approaches to circumvent our lower bounds. The formal proofs 
of our results appear in Appendix \ref{app:proofs}. 

\section{Strongly Convex and Smooth Optimization with a Second-Order Oracle}\label{sec:single}

Before presenting our main results for finite-sum optimization problems, we consider the simpler problem of minimizing a single strongly-convex and smooth function $F$ (or equivalently, \eqref{eq:obj} when $n=1$), and prove a result which may be of independent interest. 

To formalize the setting, we follow a standard oracle model, and assume that the algorithm does not have a-priori information on the objective function $F$, except the strong-convexity parameter $\lambda$ and smoothness parameter $\mu$. Instead, it has access to an oracle, which given a point $\bw\in\Wcal$, returns values and derivatives of $F$ at $\bw$ (either $\nabla F(\bw)$ for a first-order oracle, or $\nabla F(\bw),\nabla^2 F(\bw)$ for a second-order oracle). The algorithm sequentially queries the oracle using $\bw_1,\bw_2,\ldots,\bw_{T-1}$, and returns the point $\bw_T$. Our goal is to lower bound the number of oracle calls $T$, required to ensure that $\bw_T$ is an $\epsilon$-suboptimal solution. 

Given a first-order oracle and a strongly convex and smooth objective in sufficiently high dimensions, it is well-known that the worst-case oracle complexity is $$\Omega(\sqrt{\mu/\lambda}\cdot\log(1/\epsilon))$$ \citep{YudNem83}. What if we replace this by a second-order oracle, which returns both $\nabla^2 F(\bw)$ on top of $F(\bw),\nabla F(\bw)$?

Perhaps unexpectedly, it turns out that this additional information does not substantially improve the worst-case oracle complexity bound, as evidenced by the following theorem:
\begin{theorem}\label{thm:single}
	For any $\mu,\lambda$ such that $\mu>8\lambda>0$, any $\epsilon\in (0,1)$, and any deterministic algorithm, there exists a $\mu$-smooth, $\lambda$ strongly-convex function $F$ on $\reals^d$ (for $d=\tilde{\Ocal}(\sqrt{\mu/\lambda})$, hiding factors logarithmic in $\mu,\lambda,\epsilon$), such that the number of calls $T$ to a second-order oracle, required to ensure that $F(\bw_T)-\min_{\bw\in\reals^d} F(\bw)~\leq \epsilon\cdot \left(F(\mathbf{0})-\min_{\bw\in\reals^d} F(\bw)\right)$, must be at least
	\[
	c\left(\sqrt{\frac{\mu}{8\lambda}}-1\right)\cdot \log\left(\frac{(\lambda/\mu)^{3/2}}{c'\epsilon}\right),
	\]
	where $c,c'$ are positive universal constants.
\end{theorem}
For sufficiently large $\frac{\mu}{\lambda}$ and small $\epsilon$, this complexity lower bound is $\Omega\left(\sqrt{\frac{\mu}{\lambda}}\cdot \log\left(\frac{1}{\epsilon}\right)\right)$, which matches existing lower and upper bounds for optimizing strongly-convex and smooth functions using first-order methods. As mentioned earlier, the observation that such first-order oracle bounds can be extended to higher-order oracles is also briefly mentioned (without proof) in \citet[Section 7.2.6]{YudNem83}. Also, the theorem considers deterministic algorithms (which includes standard second-order methods, such as the Newton method), but otherwise makes no assumption on the algorithm. Generalizing this result to randomized algorithms should be quite doable, based on the techniques developed in \citet{woodworth2016tight}. We leave a formal derivation to future work. 

Although this result may seem surprising at first, it has a simple explanation: In order for Hessian information, which is local in nature, to be useful, there should be some regularity constraint on the Hessian, which ensures that it cannot change arbitrarily quickly as we move around the domain. A typical choice for a constraint of this kind is Lipschitz continuity which dictates that 
 $$\norm{\nabla^2 F(\bw)-\nabla^2 F(\bw')}\leq L\norm{\bw-\bw'},$$ for some constant $L$. Indeed, the construction relies on a function which does not have Lipschitz Hessians: It is based on a standard lower bound construction for first-order oracles, but the function is locally ``flattened'' in certain directions around points which are to be queried by the algorithm. This is done in such a way, that the Hessian observed by the algorithm does not provide more information than the gradient, and cannot be used to improve the algorithm's performance.

\section{Second-Order Oracle Complexity Bounds for Finite-Sum Problems}\label{sec:main}

We now turn to study finite-sum problems of the form given in \eqref{eq:obj}, and provide lower bounds on the number of oracle calls required to solve them, assuming a second-order oracle. To adapt the setting to a finite-sum problem, we assume that the second-order oracle is given both a point $\bw$ and an index $i\in \{1,\ldots,n\}$, and returns $\{f_i(\bw),\nabla f_i(\bw),\nabla^2 f_i(\bw)\}$. The algorithm iteratively produces and queries the oracle with point-index pairs $\{(\bw_t,i_t)\}_{t=1}^{T}$, with the goal of making the suboptimality (or expected suboptimality, if the algorithm is randomized) smaller than $\epsilon$ using a minimal number of oracle calls $T$.

In fact, the lower bound construction we use is such that each function $f_i$ is quadratic. Unlike the construction of the previous section, such functions have a constant (and hence trivially Lipschitz) Hessian. Moreover, since any $p$-order derivative of a quadratic for $p>2$ is zero, this means that our lower bounds automatically hold even if the oracle provides $p$-th order derivatives at any $\bw$, for arbitrarily large $p$. 

However, in order to provide a lower bound using quadratic functions, it is 
necessary to pose additional assumptions on the structure of the algorithm 
(unlike \thmref{thm:single} which is purely information-based). To see why, 
note that without computational constraints, the algorithm can simply query the 
Hessians and gradients of each $f_i(\bw)$ at $\bw=\mathbf{0}$, take the average 
to get $\nabla F(\mathbf{0})=\frac{1}{n}\sum_{i=1}^{n}\nabla f_i(\mathbf{0})$ 
and $\nabla^2 F(\mathbf{0})=\frac{1}{n}\sum_{i=1}^{n}\nabla^2 f_i(\mathbf{0})$, 
and return the exact optimum, which for quadratics equals $-\nabla^2 
F(\mathbf{0})^{-1}\nabla F(\mathbf{0})$. Therefore, with second-order 
information, the best possible information-based lower bound for quadratics is 
no better than $\Omega(n)$. This is not a satisfying bound, since in order to 
attain it we need to invert the (possibly high-rank) $d\times d$ matrix 
$\nabla^2 F(\mathbf{0})$. Therefore, if we are interested in bounds for 
computationally-efficient algorithms, we need to forbid such operations.	

Specifically, we will consider two algorithmic assumptions, which are stated 
below (their applicability to existing algorithms is discussed in the next 
section). The first assumption constrains the algorithm to query and 
return points $\bw$ which are computable using linear-algebraic manipulations 
of previous points, gradients and Hessians. Moreover, these manipulations can 
only depend on (at most) the last $\lfloor n/2 \rfloor$ Hessians returned by 
the oracle. As discussed previously, this assumption is necessary to prevent 
the algorithm from computing and inverting the full Hessian of $F$, which is 
computationally prohibitive. Formally, the assumption is the following:
\begin{assumption}[Linear-Algebraic Computations]\label{assump:linalg}
	$\bw_t$ belongs to the set $\Wcal_t\subseteq \reals^d$, defined recursively 
	as follows: $\Wcal_{1}=\{\mathbf{0}\}$, and $\Wcal_{t+1}$ is the closure of 
	the set of vectors derived from $\Wcal_{t}\cup \{\nabla f_{i_t}(\bw_t)\}$ 
	by a finite number of operations of the following form:
	\begin{itemize}
		\item $\bw,\bw' \rightarrow \alpha \bw+\alpha' \bw'$, where 
		$\alpha,\alpha'$ are arbitrary scalars.
		\item $\bw\rightarrow H\bw$, where $H$ is any $d\times d$ matrix which 
		has the same block-diagonal structure as
		\begin{equation}\label{eq:hess}
		\sum_{\tau=\max\{1,t-\lfloor n/2\rfloor+1\}}^{t}\alpha_\tau \nabla^2 
		f_{i_{\tau}}(\bw_{\tau}),
		\end{equation}
		for some arbitrary $\{\alpha_\tau\}$.
		%+D\right) \bw$.
		%\item $\bw\rightarrow \left(\sum_{\tau=\max\{1,t-\lfloor 
		%n/2\rfloor+1\}}^{t}\alpha_\tau \nabla^2 
		%f_{i_{\tau}}(\bw_{\tau})+D\right)^{-1} \bw$, assuming the matrix 
		%inverse exists.
	\end{itemize}
\end{assumption}

The first bullet allows to take arbitrary linear combinations of previous 
points and gradients, and already covers standard first-order methods and their 
variants. As to the second bullet, by ``same block-diagonal structure'', we 
mean that if the matrix in \eqref{eq:hess} can be decomposed to $r$ diagonal 
blocks of size $d_1,\ldots,d_r$ in order, then $H$ can also be decomposed into 
$r$ blocks of size $d_1,\ldots,d_r$ in order (note that this does not exclude 
the possibility that each such block is composed of additional sub-blocks). To 
give a few examples, if we let $H_t$ be the matrix in \eqref{eq:hess}, then we 
may have:
\begin{itemize}
	\item $H=H_t$,
	\item $H=H_t^{-1}$ if $H_t$ is invertible, or its pseudoinverse,
	\item $H=(H_t+D)^{-1}$ (where $D$ is some arbitrary diagonal matrix, possibly acting as a 
regularizer),
\item  $H$ is a truncated SVD decomposition of $H_t$ (or again, 
$H_t+D$ or $(H_t+D)^{-1}$ for some arbitrary diagonal matrix $D$) or its 
pseudoinverse.
\end{itemize}
Moreover, for quadratic functions, it is easily verified that
the assumption also allows prox operations (i.e. returning $\arg\min_{\bw} 
f_i(\bw)+\frac{\rho}{2}\norm{\bw-\bw'}^2$ for some $\rho,i$ and previously 
computed point $\bw'$). Also, 
note that the assumption places no limits on the number of such operations 
allowed between oracle calls. However, crucially, all these operations can be 
performed starting from a linear combination of at most $\lfloor n/2 \rfloor$ 
recent Hessians. As mentioned earlier, this is necessary, since if we could 
compute the average of all Hessians, then we could implement the Newton method. 
The assumption that the algorithm only ``remembers'' the last $\lfloor n/2 
\rfloor$ Hessians is also realistic, as existing computationally-efficient 
methods seek to use much fewer than $n$ individual Hessians at a time. We note 
that the choice of $\lfloor n/2 \rfloor$ 
is rather arbitrary, and can be replaced by $\alpha n$ for any constant 
$\alpha\in (0,1)$. Also, the way the assumption is 
formulated, the algorithm is assumed to be initialized at the origin 
$\mathbf{0}$. However, this is merely for simplicity, and can be replaced by 
any other fixed vector (the lower bound will hold by shifting the constructed 
``hard'' function appropriately).

The second (optional) assumption we will consider constrains the indices chosen 
by the 
algorithm to be oblivious, in the following sense:
\begin{assumption}[Index Obliviousness]\label{assump:obliv}
	The indices $i_1,i_2,\ldots$ chosen by the algorithm are independent of $f_1,\ldots,f_n$. 
\end{assumption}
To put this assumption differently, the indices may just as well be chosen 
before the algorithm begins querying the oracle. This can include, for 
instance, sampling functions $f_i$ uniformly at random from $f_1,\ldots,f_n$, 
and performing deterministic passes over $f_1,\ldots,f_n$ in order. As we will 
see later on, this assumption is not strictly necessary, and can be removed at 
the cost of a somewhat weaker result. Nevertheless, the assumption covers all 
optimal first-order algorithms, as well as most second-order methods we are 
aware of (see \secref{sec:related} for more details).
%Some second-order methods are based on non-uniform sampling schemes, which are 
%covered by the following relaxation:
%\begin{customassmp}{1a}[Partial Index Obliviousness]\label{assump:oblivrelax}
%	The indices $i_1,\ldots,i_T$ chosen by the algorithm depend on 
%$f_1,\ldots,f_n$ only through $\{\nabla^2 f_1(\bw),\ldots,\nabla^2 
%f_n(\bw)\}_{\bw\in\Wcal}$.
%\end{customassmp}
%We will show that one can use Assumption \ref{assump:oblivrelax} in lieu of 
%Assumption \ref{assump:obliv}, at the cost of a slightly more complicated 
%construction.
%Nevertheless, we will show that

With these assumptions stated, we can finally turn to present the main result of this section: 
\begin{theorem}\label{thm:main}
	For any $n>1$, any $\mu>\lambda>0$, any $\epsilon\in (0,c)$ (for some 
	universal constant $c>0$), and any (possibly randomized) algorithm 
	satisfying Assumptions \ref{assump:linalg} and \ref{assump:obliv}, there 
	exists $\mu$-smooth, $\lambda$-strongly convex quadratic functions 
	$f_1,\ldots,f_n:\reals^d\rightarrow \reals$ (for 
	$d=\tilde{\Ocal}(1+\sqrt{\mu/\lambda n})$, hiding factors logarithmic in 
	$n,\mu,\lambda,\epsilon$), such that the number of calls $T$ to a 
	second-order oracle, so that $$\E\left[F(\bw_T)-\min_{\bw\in\reals^d} 
	F(\bw)\right]\leq \epsilon\cdot\left(F(\mathbf{0})-\min_{\bw\in\reals^d} 
	F(\bw)\right),$$ must be at least
	\[
	\Omega\left(n+\sqrt{\frac{n\mu}{\lambda }}\cdot\log\left(\frac{(\lambda/\mu)^{3/2}\sqrt{n}}{\epsilon}\right)\right).
	\]
\end{theorem}
Comparing this with the (tight) first-order oracle complexity bounds discussed in the introduction, we see that the lower bound is the same up to log-factors, despite the availability of second-order information. In particular, the lower bound exhibits none of the favorable properties associated with full second-order methods, which can compute and invert Hessians of $F$: Whereas the full Newton method can attain $\Ocal(\log\log(1/\epsilon))$ rates, and be independent of $\mu,\lambda$ if $F$ satisfies a self-concordance property \citep{boyd2004convex}, here we only get a linear $\Ocal(\log(1/\epsilon))$ rate, and there is a strong dependence on $\mu,\lambda$, even though the function is quadratic and hence self-concordant.

The proof of the theorem is based on a randomized construction, which can be sketched as follows: We choose indices $j_1,\ldots,j_{d-1}\in \{1,\ldots,n\}$ independently and uniformly at random, and define
\begin{align*}
f_i(\bw) = ~a&\cdot w_1^2+\hat{a}\cdot \sum_{l=1}^{d-1}\mathbf{1}_{j_l=i}(w_l-w_{l+1})^2\\
&+\bar{a}\cdot w_d^2-\tilde{a}\cdot  w_1+\frac{\lambda}{2}\norm{\bw}^2,
\end{align*}
where $\mathbf{1}_A$ is the indicator function of the event $A$, and  $a,\hat{a},\bar{a},\tilde{a}$ are parameters chosen based on $\lambda,\mu,n$. 
The average function $F(\bw)=\frac{1}{n}\sum_{i=1}^{n}f_i(\bw)$ equals
\[
F(\bw) = a\cdot w_1^2+\frac{\hat{a}}{n}\cdot \sum_{l=1}^{d-1}(w_l-w_{l+1})^2+\bar{a}\cdot w_d^2-\tilde{a}\cdot  w_1+\frac{\lambda}{2}\norm{\bw}^2.
\]
By setting the parameters appropriately, it can be shown that $F$ is $\lambda$-strongly convex and each $f_i$ is $\mu$-smooth. Moreover, the optimum of $F$ has the form $(q,q^2,q^3,\ldots,q^d)$ for 
$$q=\frac{\sqrt\kappa-1}{\sqrt\kappa+1},$$ where 
\begin{align} \label{eq:kappa}
\kappa=\frac{\frac\mu\lambda-1}n+1
\end{align}
is the so-called condition number of $F$. The proof is based on arguing that after $T$ oracle calls, the points computable by any algorithm satisfying Assumptions \ref{assump:obliv} and \ref{assump:linalg} must have $0$ values at all coordinates larger than some $l_T$, hence the squared distance of $\bw_T$ from the optimum must be at least $\sum_{i=l_T+1}^{d} q^{2i}$, which leads to our lower bound. Thus, the proof revolves around upper bounding $l_T$. We note that a similar construction of $F$ was used in some previous first-order lower bounds under algorithmic assumptions (e.g. \citet{nesterov2013introductory,lan2015optimal}, as well as \citet{arjevani2015communication} in a somewhat different context). The main difference is in how we construct the individual functions $f_i$, and in analyzing the effect of second-order rather than just first-order information. 

To upper bound $l_T$, we let $l_t$ (where $t=1,\ldots,T$) be the largest non-zero coordinate in $\bw_t$, and track how $l_t$ increases with $t$. The key insight is that if $\bw_1,\ldots,\bw_{t-1}$ are zero beyond some coordinate $l$, then any linear combinations of them, \emph{as well as multiplying them by matrices based on second-order information}, as specified in Assumption \ref{assump:linalg}, will still result in vectors with zeros beyond coordinate $l$. The only way to ``advance'' and increase the set of non-zero coordinates is by happening to query the function $f_{j_l}$. However, since the indices of the queried functions are chosen obliviously, whereas each $j_l$ is chosen uniformly at random, the probability of this happening is quite small, of order $1/n$. Moreover, we show that even if this event occurs, we are unlikely to ``advance'' by more than $\Ocal(1)$ coordinates at a time. Thus, the algorithm essentially needs to make $\Omega(n)$ oracle calls in expectation, in order to increase the number of non-zero coordinates by $\Ocal(1)$. It can be shown that the number of coordinates needed to get an $\epsilon$-optimal solution is $\tilde{\Omega}(\sqrt{\mu/n\lambda}\cdot \log(1/\epsilon))$ (hiding some log-factors). Therefore, the total number of oracle calls is about $n$ times larger, namely $\tilde{\Omega}(\sqrt{n\mu/\lambda}\cdot \log(1/\epsilon))$. To complete the theorem, we also provide a simple and separate $\Omega(n)$ lower bound, which holds since each oracle call gives us information on just one of the $n$ individual functions $f_1,\ldots,f_n$, and we need some information on most of them in order to get a close-to-optimal solution.

When considering non-oblivious (i.e., adaptive) algorithms, the construction 
used in \thmref{thm:main} fails as soon as the algorithm obtains the Hessians 
of all the individual functions (potentially, after $n$ oracle queries). 
Indeed, knowing the Hessians of $f_i$, one can devise an index-schedule which 
gains at least one coordinate at every iteration (by querying the function 
which holds the desired $2\times2$ block), as opposed to $O(1/n)$ on average in 
the oblivious case. Nevertheless, as mentioned before, we can still provide a 
result similar to \thmref{thm:main} even if the indices are chosen adaptively, 
at the cost of a much larger dimension:
\begin{theorem}\label{thm:mainrelaxed}
	\thmref{thm:main} still holds if one omits Assumption \ref{assump:obliv}, 
	and with probability $1$ rather than in expectation, at 
	the cost of requiring an exponentially larger dimensionality of $d = 
	n^{\tilde{\Ocal}\left({1+\sqrt{\mu/\lambda n}}\right)}$.
\end{theorem}
The proof is rather straightforward: Making the dependence on the random 
indices $j_1,\ldots,j_{d-1}$ explicit, the quadratic construction used in the 
previous theorem can be written as 
\begin{align*}
F^{j_1,\ldots,j_{d-1}}&(\bw)=\frac{1}{n}\sum_{i=1}^{n}f^{j_1,\ldots,j_{d-1}}_i(\bw)\\
&= \frac{1}{n}\sum_{i=1}^{n}\bw^\top A^{j_1,\ldots,j_{d-1}}_i\bw-\tilde{a}\inner{\be_1,\bw}+\frac{\lambda}{2}\norm{\bw}^2
\end{align*}
for some $d\times d$ matrix $A^{j_1,\ldots,j_{d-1}}_i$ dependent on 
$j_1,\ldots,j_{d-1}$, and a fixed parameter $\tilde{a}$. Now, we create $n$ 
huge block-diagonal matrices 
$A_1,\ldots,A_n$, where each $A_i$ contains $A^{j_1,\ldots,j_{d-1}}_i$ for each 
of the $n^{d-1}$ possible choices of $j_1,\ldots,j_{d-1}$ along its diagonal 
(in some canonical order), and one huge vector 
$$\be = \be_1+\be_{d+1}+\ldots+\be_{(n^{d-1}-1)d+1}.$$ We then let
\begin{align*}
F(\bw) &= \frac{1}{n}\sum_{i=1}^{n}f_i(\bw)\\ &= \frac{1}{n}\sum_{i=1}^{n}\bw^\top 
A_i \bw -\tilde{a}\inner{\be,\bw}+\frac{\lambda}{2}\norm{\bw}^2.	
\end{align*}
This function essentially combines all $n^{d-1}$ problems 
$F^{j_1,\ldots,j_{d-1}}$ simultaneously, where each $F^{j_1,\ldots,j_{d-1}}$ 
is embedded in a disjoint set of coordinates. Due to the block-diagonal 
structure of each $A_i$, this function inherits the strong-convexity and 
smoothness properties of the original construction. Moreover, to optimize this 
function, the algorithm needs to ``solve'' all $n^{d-1}$ problems 
simultaneously, using the same choice of indices $i_1,i_2,\ldots$. Using a 
combinatorial argument which parallels the probabilistic argument in the proof 
of \thmref{thm:main}, we can show that no matter how these indices are chosen, 
the average number of non-zero coordinates of the iterates cannot grow too 
rapidly, and lead to the same bound as in \thmref{thm:main}. Since the 
construction is deterministic, and applies no matter how the indices are 
chosen, the lower bound holds deterministically, rather than in expectation as 
in \thmref{thm:main}.

Lastly, it is useful to consider how the bounds stated in \thmref{thm:main} and 
\thmref{thm:mainrelaxed} differ when the dimension $d$ is fixed and finite. 
Inspecting the proofs of both theorems reveals that in both cases the 
suboptimality, as a function of the iteration number $T$, has a linear 
convegence rate bounded from below by
\begin{align} \label{ineq:conv_rate}
\E\left[\frac{F(\bw_T)-F(\bw^\star)}{F(\mathbf{0})-F(\bw^\star)}\right]\geq \Omega(1)\left(\frac{\sqrt{\kappa}-1}{\sqrt{\kappa}+1}\right)^{\Ocal\left(\frac{T}{n}\right)}
\end{align}
(where $\kappa$ is as defined in \eqref{eq:kappa}, and $\Omega(1)$ hides 
dependencies on the problem parameters, but is independent of $T$).
However, whereas the bound established in \thmref{thm:main} is valid for 
$\Ocal(d)$ number of iterations, \thmref{thm:mainrelaxed} applies to a much 
restricted range of roughly $\log(d)/\log(n)$ iterations. This indicate 
that adaptive optimization algorithms might be able to gain a super-linear 
convergence rate after a significantly smaller number of iterations in 
comparison to oblivious algorithms (see \cite{arjevani2016dimension} for a 
similar discussion regarding first-order methods). That being said, trading 
obliviousness for adaptivity may increase the per-iteration cost and reduce 
numerical stability. 

%depending on the 
%choice of $r$, the function is equivalent to the construction involving some 
%$j_1,\ldots,j_{d-1}$ on some subset of $d$ coordinates (and along the other 
%coordinates, the optimum is simply $0$). The only main difference is that the 
%quadratic terms $A_i$ are no longer random -- all the randomness is ``pushed'' 
%to the linear term $\be_r$, and the random $r$ encodes which choices of 
%$j_1,\ldots,j_{d-1}$ are used. Thus, even if the indices $i_1,i_2,\ldots$ used 
%by the algorithm depends on the Hessians $A_i$, they are independent of $r$, 
%and hence independent of $j_1,\ldots,j_{d-1}$ corresponding to the 
%construction. Therefore, the analysis of \thmref{thm:main} (which assumes such 
%independence via Assumption \ref{assump:obliv}) holds verbatim, and the result 
%follows.

\section{Comparison to Existing Approaches}\label{sec:related}

As discussed in the introduction, there has been a recent burst of activity involving second-order methods for solving finite-sum problems, relying on Hessians of individual functions $f_i$. In this section, we review the main algorithmic approaches and compare them to our results. The bottom line is that most existing approaches satisfy the assumptions stated in \secref{sec:main}, and therefore our lower bounds will apply, at least in a worst-case sense. A possible exception to this is the Newton sketch algorithm \citep{pilanci2015newton}, which relies on random projections, but on the flip side is computationally expensive. 

Turning to the details, existing approaches are based on taking the standard Newton iteration for such problems, 
\begin{align*}
\bw&_{t+1}=\bw_t-\alpha_t \left(\nabla^2 F(\bw_t)\right)^{-1}\nabla F(\bw_t)\\
&=\bw_t-\alpha_t \left(\frac{1}{n}\sum_{i=1}^{n}\nabla^2 
f_i(\bw_t)\right)^{-1}\left(\frac{1}{n}\sum_{i=1}^{n}\nabla f_i(\bw_t)\right),
\end{align*}
and replacing the inverse Hessian term $\left(\frac{1}{n}\sum_{i=1}^{n}\nabla^2 f_i(\bw)\right)^{-1}$ (and sometimes the vector term $\frac{1}{n}\sum_{i=1}^{n}\nabla f_i(\bw)$ as well) by some approximation which is computationally cheaper to compute. One standard and well-known approach is to use only gradient information to construct such an approximation, leading to the family of quasi-Newton methods \citep{nocedal2006numerical}. However, as they rely on first-order rather than second-order information, they are orthogonal to the topic of our work, and are already covered by existing  complexity lower bounds for first-order oracles.

Turning to consider Hessian approximation techniques using second-order information, perhaps the simplest and most intuitive approach is sampling: Since the Hessian equals the average of many individual Hessians,  
$$\nabla^2 F(\bw)=\frac{1}{n}\sum_{i=1}^{n}\nabla^2 f_i(\bw),$$ we can approximate it by taking a sample $S$ of indices in $\{1,\ldots,n\}$ uniformly at random, compute the Hessians of the corresponding individual functions, and use the approximation
\[
\nabla^2 F(\bw)\approx\frac{1}{|S|}\sum_{i\in S}\nabla^2 f_i(\bw).
\]
If $|S|$ is large enough, then by concentration of measure arguments, this sample average should be close to the actual Hessian $\nabla^2 F(\bw)$. On the other hand, if $|S|$ is not too large, then the resulting matrix is easier to invert (e.g. because it has a rank of only $\Ocal(|S|)$, if each individual Hessian has rank $\Ocal(1)$, as in the case of linear predictors). Thus, one can hope that the right sample size will lead to computational savings. There have been several rigorous studies of such ``subsampled Newton'' methods, such as \citet{erdogdu2015convergence,roosta2016sub,roosta2016sub2,bollapragada2016exact} and references therein. However, our lower bound in  \thmref{thm:main} holds for such an approach, since it satisfies both Assumption \ref{assump:obliv} and \ref{assump:linalg}. As expected, the existing worst-case complexity upper bounds are no better than our lower bound.

\cite{xu2016sub} recently proposed a subsampled Newton method, together with 
\emph{non-uniform} sampling, which assigns more weight to individual functions 
which are deemed more ``important''. This is measured via properties of the 
Hessians of the functions, such as their norms or via leverage scores. This 
approach breaks Assumption \ref{assump:obliv}, as the sampled indices are now 
chosen in a way dependent on the individual functions. However, our lower bound 
in \thmref{thm:mainrelaxed}, which does not require this assumption, still 
applies to such a method. 

A variant of the subsampled Newton approach, studied in \citet{erdogdu2015convergence}, uses a low-rank approximation of the sample Hessian (attained by truncated SVD), in lieu of the sample Hessian itself. However, this still falls in the framework of Assumption \ref{assump:linalg}, and our lower bound still applies.

A different approach to approximate the full Hessian is via randomized 
sketching techniques, which replace the Hessian $\nabla^2 F(\bw)$ by a low-rank 
approximation of the form $$(\nabla^2 F(\bw))^{1/2}SS^\top (\nabla^2 
F(\bw))^{1/2},$$ where $S\in \reals^{d\times m},m\ll d$ is a random sketching 
matrix, and $ (\nabla^2 F(\bw))^{1/2}$ is the matrix square root of $\nabla^2 
F(\bw)$. This approach forms the basis of the Newton sketch algorithm proposed 
in \citet{pilanci2015newton}. This approach currently escapes our lower bound, 
since it violates Assumption \ref{assump:linalg}. That being said, this 
approach is inherently expensive in terms of computational resources, as it 
requires us to compute the square root of the full Hessian matrix. Even under 
favorable conditions, this requires us to perform a full pass over all 
functions $f_1,\ldots,f_n$ at every iteration. Moreover, existing iteration 
complexity upper bounds have a strong dependence on both $\mu/\lambda$ as well 
as the dimension $d$, and are considerably worse than the lower bound of 
\thmref{thm:main}. Therefore, we conjecture that this approach cannot lead to 
better worst-case results.

\citet{agarwal2016second} develop another line of stochastic second-order methods, which are based on the observation that the Newton step $(\nabla^2 F(\bw))^{-1}\nabla F(\bw)$ is the solution of the system of linear equations $$\nabla^2 F(\bw)\bx=\nabla F(\bw).$$ Thus, one can reduce the optimization problem to solving this system as efficiently as possible. The basic variant of their algorithm (denoted as LiSSA) relies on operations of the form $$\bw\mapsto (I-\nabla^2 f_i(\bw))\bw$$ (for $i$ sampled uniformly at random), as well as linear combinations of such vectors, which satisfy our assumptions. A second variant, LiSSA-Quad, re-phrases this linear system as the finite-sum optimization problem
\begin{align*}
\min_{\bx}~ \bx^\top &\nabla^2 F(\bw)\bx+\nabla F(\bw)^\top \bx\\
&= \frac{1}{n}\sum_{i=1}^{n}\bx^\top \nabla^2 f_i(\bw)\bx+\nabla f_i(\bw)^\top \bx,
\end{align*}
and uses some first-order method for finite-sum problems in order to solve it. Since individual gradients of this objective are of the form $\nabla^2 f_i(\bw)\bx+\nabla f_i(\bw)$, and most state-of-the-art first-order methods pick indices $i$ obliviously, this approach also satisfies our assumptions, and our lower bounds apply. Yet another proposed algorithm, LiSSA-Sample, is based on replacing the optimization problem above by
\begin{equation}\label{eq:objBtrans}
\min_{\bx}~ \bx^\top \nabla^2 F(\bw)B^{-1}\bx+\nabla F(\bw)^\top \bx, 
\end{equation}
where $B$ is some invertible matrix, solving it (with the optimum being equal 
to $B(\nabla^2 F(\bw))^{-1}\nabla F(\bw)$), and multiplying the solution by 
$B^{-1}$ to recover the solution $(\nabla^2 F(\bw))^{-1}\nabla F(\bw)$ to the 
original problem. In order to get computational savings, $B$ is chosen to be a 
linear combination of $\Ocal(d\log(d))$ sampled individual hessians $\nabla^2 
f_i(\bw)$, where it is assumed that $d\log(d)\ll n$, and the sampling and 
weighting is carefully chosen (based on the Hessians) so that 
\eqref{eq:objBtrans} has strong convexity and smoothness parameters within a 
constant of each other. As a result, \eqref{eq:objBtrans} can be solved quickly 
using standard gradient descent, taking steps along the gradient, which equals 
$\nabla^2 F(\bw)B^{-1}\bx+\nabla F(\bw)$ at any point $\bx$. This gradient is 
again computable under \ref{assump:linalg} (using $\Ocal(n)$ oracle calls), 
since $B$ is a linear combination of $d\log(d)\ll n$ sampled individual 
Hessians. Thus, our lower 
bound (in the form of \thmref{thm:mainrelaxed}) still applies to such methods. 

That being said, it is important to note that the complexity upper bound 
attained in \citet{agarwal2016second} for LiSSA-Sample is on the order of 
$$\tilde{\Ocal}((n+\sqrt{d\mu/\lambda})\cdot\text{polylog}(1/\epsilon))$$ 
(at least asymptotically as $\epsilon\rightarrow 0$), which can be better than 
our lower 
bound if $d\ll n$. There is no contradiction, since the lower bound in 
\thmref{thm:mainrelaxed}  only applies for a dimension $d$ much larger than 
$n$. Interestingly, our results also indicate that an adaptive index
sampling scheme is \emph{necessary} to get this kind of improved performance 
when $d\ll n$: Otherwise, it could violate \thmref{thm:main}, which establishes 
a lower bound of
$\tilde{\Ocal}(n+\sqrt{n\mu/\lambda})$ even if the dimension is quite moderate 
($d=\tilde{\Ocal}(1+\sqrt{\mu/\lambda n})$, which is $\ll n$ under the mild 
assumption that $\mu/\lambda\ll n^3$). 

The observation that an adaptive scheme (breaking assumption 
\ref{assump:obliv}) can help performance when $d\ll n$ is also  seen in the 
lower bound construction used to prove \thmref{thm:main}: If $\mu,\lambda,n$ 
are such that the required dimension $d$ is $\ll n$, then it means that only 
the functions $f_{j_1},\ldots,f_{j_{d-1}}$, which are a small fraction of all 
$n$ individual functions, are informative and help us reduce the objective 
value. Thus, sampling these functions in an adaptive manner is imperative to 
get better complexity than the bound in \thmref{thm:main}. Based on the fact 
that only at most $d-1$ out of $n$ functions are relevant in the construction, 
we conjecture that the possible improvement in the worst-case oracle complexity 
of such 
schemes may amount to replacing dependencies on $n$ with dependencies on $d$, 
which is indeed the type of improvement attained (for small enough $\epsilon$) 
in \citet{agarwal2016second}.

Finally, we note that \citet{agarwal2016second} proposes another algorithm tailored to self-concordant functions, with runtime independent of the smoothness and strong convexity parameters of the problem. However, it requires performing $\geq 1$ full Newton steps, so the runtime is prohibitive for large-scale problems (indeed, for quadratics as used in our lower bounds, even a single Newton step suffices to compute an exact solution).

\section{Summary and Discussion}\label{sec:summary}

In this paper, we studied the oracle complexity for optimization problems, assuming availability of a second-order oracle. This is in contrast to most existing oracle complexity results, which focus on a first-order oracle. First, we formally proved that in the standard setting of strongly-convex and smooth optimization problems, second-order information does not significantly improve the oracle complexity, and further assumptions (i.e. Lipschitzness of the Hessians) are in fact necessary. We then presented our main lower bounds, which show that for finite-sum problems with a second-order oracle, under some reasonable algorithmic assumptions, the resulting oracle complexity is -- again -- not significantly better than what can be obtained using a first-order oracle. Moreover, this is shown using quadratic functions, which have $0$ derivatives of order larger than $2$. Hence, our lower bounds apply even if we have access to an oracle returning derivatives of order $p$ for all $p\geq 0$, and the function is smooth to any order. In \secref{sec:related}, we studied how our framework and lower bounds are applicable to most existing approaches.

Although this conclusion may appear very pessimistic, they are actually useful in pinpointing potential assumptions and approaches which may circumvent these lower bounds. In particular:
\begin{itemize}
	\item Our lower bound for algorithms employing adaptive index sampling 
	schemes (\thmref{thm:mainrelaxed}) only hold when the dimension $d$ is very 
	large. This leaves open the possibility of better (non index-oblivious) 
	algorithms when $d$ is moderate, as was recently demonstrated in the 
	context of the LiSSA-Sample algorithm of \citet{agarwal2016second} (at 
	least for small enough $\epsilon$). As discussed in the previous section, 
	we conjecture that the possible improvement in the worst-case oracle 
	complexity of 
	such schemes may amount to replacing dependencies on $n$ with dependencies 
	on $d$.
	\item It might be possible to construct algorithms breaking Assumption 
	\ref{assump:linalg}, e.g. by using operations which are not 
	linear-algebraic. That being said, we currently conjecture that this 
	assumptions 
	can be 
	significantly relaxed, and similar results would hold for any algorithm 
	which has ``significantly'' cheaper iterations (in terms of runtime) 
	compared to the Newton method. 
	\item Our lower bounds are worst-case over smooth and strongly-convex individual functions $f_i$. It could be that by assuming more structure, better bounds can be obtained. For example, as discussed in the introduction, an important special case is when $f_i(\bw)=\ell_i(\bx_i^\top \bw)$ for some scalar function $\ell_i$ and vector $\bx_i$. Our construction in \thmref{thm:main} does not quite fit this structure, although it is easy to show that we still get functions of the form $f_i(\bw)=\ell_i(X_i^\top\bw)$, where $X_i$ has $\Ocal(1+d/n)=\tilde{\Ocal}(1+\sqrt{\mu/\lambda n^3})$ rows in expectation, which is $\tilde{\Ocal}(1)$ under a broad parameter regime. We believe that the difference between $\tilde{\Ocal}(1)$ rows and $1$ row is not significant in terms of the attainable oracle complexity, but we may be wrong. Another possibility is to provide results depending on more delicate spectral properties of the function, beyond its strong convexity and smoothness, which may lead to better results and algorithms under favorable assumptions.
	\item Our lower bounds in \secref{sec:main}, which establish a linear 
convergence rate (logarithmic dependence on $\log(1/\epsilon)$), are 
non-trivial only if the optimization error $\epsilon$ is sufficiently small. 
This does not preclude the possibility of attaining better initial performance 
when $\epsilon$ is relatively large. 
%Indeed, some of the recent work mentioned 
%previously (e.g. 
%\cite{bollapragada2016exact,xu2016sub,roosta2016sub,pilanci2015newton}) show 
%that various algorithms can have an initial super-linear convergence rate, 
%which slows down to linear as one gets closer to the optimum. 
\end{itemize}

In any case, we believe our work lays the foundation for a more comprehensive 
study of the complexity of efficient second-order methods, for finite-sum and 
related optimization and learning problems.  

\subsubsection*{Acknowledgments}
This research is supported in part by an FP7 Marie Curie CIG grant and an Israel Science Foundation grant 425/13.

%                                      IMPORTANT!! add for camera ready!
%\subsubsection*{Acknowledgments}
%This research is supported in part by an FP7 Marie Curie CIG grant and an Israel Science Foundation grant 425/13.

\bibliographystyle{plainnat}
\bibliography{mybib}

\appendix
\onecolumn

\section{Proofs}\label{app:proofs}

\subsection{Auxiliary Lemmas}
The following lemma was essentially proven in \cite{lan2015optimal,nesterov2013introductory}, but we provide a proof for completeness:

\begin{lemma}\label{lem:tridig}
	Fix $\alpha,\beta \ge 0$, and consider the following function on $\reals^d$:
	\[
	F(\bw) = \frac{\alpha}{8}\left(w_1^2+\sum_{i=1}^{d-1}(w_i-w_{i+1})^2+(a_{\tilde{\kappa}}-1)w_d^2-w_1\right)+\frac{\beta}{2}\norm{\bw}^2,
	\]
and $a_{\tilde{\kappa}} = \frac{\sqrt{\tilde{\kappa}}+3}{\sqrt{\tilde{\kappa}}+1}$ where $\tilde{\kappa} =\frac{\alpha+\beta}{\beta}$ is the condition number of $F$. Then $F$ is $\beta$ strongly convex, $(\alpha+\beta)$-smooth, and has a unique minimum at $(q,q^2,q^3,\ldots,q^d)$ where $q=\frac{\sqrt{\tilde{\kappa}}-1}{\sqrt{\tilde{\kappa}}+1}$.
\end{lemma}
\begin{proof}
	The function is equivalent to 
	\[
	F(\bw) = \frac{\alpha}{8}\left(\bw^\top A \bw-w_1\right)+\frac{\beta}{2}\norm{\bw}^2,
	\]
	where
	\[
	A = \left(\begin{matrix}
	2 & -1 &  &  &  \\ 
	-1 & 2 & -1 &  &  \\ 
	& -1 & \ddots & \ddots &  \\ 
	&  & \ddots & 2 & -1 \\ 
	&  &  & -1 & a_{\tilde{\kappa}}
	\end{matrix} \right).
	\]
	Since $A$ is symmetric, all its eigenvalues are real. Therefore, by Gershgorin circle theorem and the fact that $a_{\tilde{\kappa}}\in[1,2]$ (since $\tilde{\kappa}\ge1$), we have that all the eigenvalues of $A$ lie in $[0,4]$. Thus, the eigenvalues of $\nabla^2 F=(\alpha/4) A + \beta I$ lie in $[\beta,\alpha+\beta]$, implying that $F$ is $\beta$-strongly convex and $(\alpha+\beta)$-smooth.
	
	It remains to compute the optimum of $F$. By differentiating $F$ and setting to zero, we get that the optimum $\bw$ must satisfy the following set of equations:
	\begin{align*}
	&w_2-2\cdot\frac{\tilde{\kappa}+1}{\tilde{\kappa}-1}\cdot w_1+1=0\\
	&w_{i+1}-2\cdot\frac{\tilde{\kappa}+1}{\tilde{\kappa}-1}\cdot w_i+w_{i-1}=0~~~~\forall~i=2,\ldots,d-1\\
	&\left(a_{\tilde{\kappa}}+\frac{4}{\tilde{\kappa}-1}\right)w_{d}-w_{d-1}=0.
	\end{align*}
	It is easily verified that this is satisfied by the vector $(q,q^2,q^3,\ldots,q^d)$, where $q=\frac{\sqrt{\tilde{\kappa}}-1}{\sqrt{\tilde{\kappa}}+1}$. Since $F$ is strongly convex, this stationary point must be the unique global optimum of $F$.
\end{proof}

\begin{lemma}\label{lem:gq_to_ell}
	For some $q\in (0,1)$ and positive $d$, define \[g(z) = \begin{cases}q^{2(z+1)}& z< d\\0 & z\geq d\end{cases}~.\] Let $l$ be a non-negative random variable, and suppose $d\geq 2\E[l]$. Then $\E[g(l)] \geq \frac{1}{2}q^{2\E[l]+2}$.
\end{lemma}
\begin{proof}
	Since $q\in (0,1)$, the function $z\mapsto q^z$ is convex for non-negative $z$ and monotonically decreasing. Therefore, by definition of $g$ and Jensen's inequality, we have
	\begin{align*}
	\E[g(l)] &= \Pr(l<d)\cdot\E[q^{2(l+1)}|l<d]+\Pr(l\geq d)\cdot 0
	\geq \Pr(l<d)\cdot q^{\E[2(l+1)]}.
	\end{align*}
	Using Markov's inequality to derive $\Pr(l<d) = 1-\Pr(l\geq d) \geq 1-\frac{\E[l]}{d} \geq \frac{1}{2}$, concludes the proof.
\end{proof}

\subsection{Proof of \thmref{thm:single}}

The proof is inspired by a technique introduced in \cite{woodworth2016tight} for analyzing randomized first-order methods, in which a quadratic function is ``locally flattened'' in order to make first-order (gradient) information non-informative. We use a similar technique to make \emph{second-order} (Hessian) information non-informative, hence preventing second-order methods from having an advantage over first-order methods.

Given a (deterministic) algorithm and a bound $T$ on the number of oracle calls, we construct the function $F$ in the following manner. We first choose some dimension $d \geq 2T$.
We then define 
\[
\kappa = \frac{\mu}{8\lambda}~~~,~~~ q=\frac{\sqrt{\kappa}-1}{\sqrt{\kappa}+1},
\]
and choose $r>0$ sufficiently small so that
\[
\frac{T\mu r^2}{8\lambda}\leq 1~~~\text{and}~~~\sqrt{\frac{T\mu r^2}{16\lambda}}\leq \frac{1}{2}q^{T}.
\]
We also let $\bv_1,\ldots,\bv_{T}$ be orthonormal vectors in $\reals^d$ (to be specified later). We finally define our function as
\[
F(\bw) = H(\bw)+\frac{\lambda}{2}\norm{\bw}^2,
\]
where
\[
H(\bw) =  \frac{\lambda(\kappa-1)}{8}\left(\inner{\bv_1,\bw}^2+\sum_{i=1}^{T-1}\phi_r(\inner{\bv_i-\bv_{i+1},\bw})+(a_{\kappa}-1)\phi_r(\inner{\bv_{T},\bw})-\inner{\bv_1,\bw}\right),
\]
$a_{\kappa} = \frac{\sqrt{\kappa}+3}{\sqrt{\kappa}+1}$, and 
\[
\phi_r(z) = \begin{cases} 0 & |z|\leq r\\ 2(|z|-r)^2 & r<|z|\leq 2r \\ z^2-2r^2 & |z|>2r\end{cases}~.
\]
It is easy to show that $\phi_r$ is $4$-smooth and satisfies $0\leq z^2-\phi_r(z) \leq 2r^2$ for all $z$.

First, we establish that $F$ is indeed strongly convex and smooth as required:
\begin{lemma}\label{lem:Fstrsmooth}
	$F$ as defined above is $\lambda$-strongly convex and $\mu$-smooth.
\end{lemma}
\begin{proof}
	Since $\phi_r$ is convex, and the composition of a convex and linear function is convex, we have that $\bw\mapsto \phi_r(\inner{\bv_i-\bv_{i+1},\bw})$ are convex for all $i$, as well as $\bw\mapsto \inner{\bv_1,\bw}^2$ and $\bw\mapsto \phi_r(\inner{\bv_T,\bw})$. Therefore, $H(\bw)$ is convex. As a result, $F$ is $\lambda$-strongly convex due to the $\frac{\lambda}{2}\norm{\bw}^2$ term. As to smoothness, note first that $H(\bw)$ can be equivalently written as $\tilde{H}(V\bw)$, where $V$ is some orthogonal $d\times d$ matrix with the first $T$ rows equal to $\bv_1,\ldots,\bv_{T}$, and
	\[
	\tilde{H}(\bx) = \frac{\lambda(\kappa-1)}{8}\left(x_1^2+\sum_{i=1}^{T-1}\phi_r(x_i-x_{i+1})+(a_{\kappa}-1)\phi_r(x_{T})-x_1\right).
	\]
	Therefore, $\nabla^2 F(\bw) = \nabla^2 H(\bw)+\lambda I = V^\top \nabla^2 \tilde{H}(V\bw)V+\lambda I$. It is easily verified that $\nabla^2 \tilde{H}$ at any point (and in particular $V\bw$) is tridiagonal, with each element having absolute value at most $2\lambda(\kappa-1)$. Therefore, using the orthogonality of $V$ and the fact that $(a+b)^2\leq 2(a^2+b^2)$,
	\begin{align*}
	\sup_{\bx:\norm{\bx}=1}\bx^\top \nabla^2 F(\bw) \bx ~&=~ \sup_{\bx:\norm{\bx}=1}\bx^\top(V^\top \nabla^2 \tilde{H}(V\bw)V+\lambda I)\bx \\
	&=~ \sup_{\bx:\norm{\bx}=1}\bx^\top \nabla^2 \tilde{H}(V\bw)\bx+\lambda\\
	&\leq~\sup_{\bx:\norm{\bx}=1} 2\lambda(\kappa-1)\left(\sum_{i=1}^{d}x_i^2+2\sum_{i=1}^{d-1}|x_ix_{i+1}|\right)+\lambda\\
	&\leq~ \sup_{\bx:\norm{\bx}=1}2\lambda(\kappa-1)\sum_{i=1}^{d-1}(|x_i|+|x_{i+1}|)^2+\lambda\\
	&\leq~ \sup_{\bx:\norm{\bx}=1}4\lambda(\kappa-1)\sum_{i=1}^{d-1}(x_i^2+x_{i+1}^2)+\lambda\\
	&\leq~ 8\lambda(\kappa-1)+\lambda~\leq~ 8\lambda\kappa.
	\end{align*}
	Plugging in the definition of $\kappa$, this equals $\mu$. Therefore, the spectral norm of the Hessian of $F$ at any point is at most $\mu$, and therefore $F$ is $\mu$-smooth.
\end{proof}

By construction, the function $F$ also has the following key property:
\begin{lemma}\label{lem:info}
	For any $\bw\in \reals^d$ orthogonal to $\bv_t,\bv_{t+1},\ldots,\bv_{T}$ (for some $t\in \{1,2,\ldots,T-1\}$), it holds that $F(\bw), \nabla F(\bw),\nabla^2 F(\bw)$ do not depend on $\bv_{t+1},\bv_{t+2},\ldots,\bv_{T}$.
\end{lemma}
\begin{proof}
	Recall that $F$ is derived from $H$ by adding a $\frac{\lambda}{2}\norm{\bw}^2$ term, which clearly does not depend on $\bv_1,\ldots,\bv_{T}$. Therefore, it is enough to prove the result for $H(\bw),\nabla H(\bw),\nabla^2 H(\bw)$. By taking the definition of $H$ and differentiating, we have that $H(\bw)$ is proportional to 
	\[
	\inner{\bv_1,\bw}^2+\sum_{i=1}^{T-1}\phi_r(\inner{\bv_i-\bv_{i+1},\bw})+(a_{\kappa}-1)\phi_r(\inner{\bv_{T},\bw})-\inner{\bv_1,\bw},
	\]
	$\nabla H(\bw)$ is proportional to
	\[
	2\inner{\bv_1,\bw}\bv_1+\sum_{i=1}^{T-1}\phi'_{r}(\inner{\bv_i-\bv_{i+1},\bw})(\bv_i-\bv_{i+1})+(a_{\kappa}-1)\phi'_{r}(\inner{\bv_{T},\bw})\bv_{T}-\bv_1,
	\]
	and $\nabla^2 H(\bw)$ is proportional to
	\[
	2\bv_1\bv_1^\top+\sum_{i=1}^{T-1}\phi''_{r}(\inner{\bv_i-\bv_{i+1},\bw})(\bv_i-\bv_{i+1})(\bv_i-\bv_{i+1})^\top+(a_{\kappa}-1)\phi''_{r}(\inner{\bv_{T},\bw})\bv_{T}\bv_{T}^\top.
	\]
	By the assumption $\inner{\bv_{t},\bw}=\inner{\bv_{t+1},\bw}=\ldots=\inner{\bv_{T},\bw}=0$, and the fact that $\phi_r(0)=\phi'_{r}(0)=\phi''_{r}(0)=0$, we have $\phi_r(\inner{\bv_i-\bv_{i+1},\bw})=\phi'_{r}(\inner{\bv_i-\bv_{i+1},\bw})=\phi''_{r}(\inner{\bv_i-\bv_{i+1},\bw})=0$ for all $i\in \{t,t+1,\ldots,T\}$, as well as $\phi_r(\inner{\bv_T,\bw})=\phi'_{r}(\inner{\bv_T,\bw})=\phi''_{r}(\inner{\bv_T,bw})=0$. Therefore, it is easily verified that the expressions above indeed do not depend on $\bv_{t+1},\ldots,\bv_{T}$. 
\end{proof}

With this lemma at hand, we now turn to describe how $\bv_1,\ldots,\bv_{T}$ are constructed:
\begin{itemize}
	\item First, we compute $\bw_1$ (which is possible since the algorithm is deterministic and $\bw_1$ is chosen before any oracle calls are made).
	\item We pick $\bv_1$ to be some unit vector orthogonal to $\bw_1$. Assuming $\bv_2,\ldots,\bv_{T}$ will also be orthogonal to $\bw_1$ (which will be ensured by the construction which follows), we have by \lemref{lem:info} that the information $F(\bw_1),\nabla F(\bw_1),\nabla^2 F(\bw_1)$ provided by the oracle to the algorithm does not depend on $\{\bv_2,\ldots,\bv_{T}\}$, and thus depends only on $\bv_1$ which was already fixed. Since the algorithm is deterministic, this fixes the next query point $\bw_2$.
	\item For $t=2,3,\ldots,T-1$, we repeat the process above: We compute $\bw_t$, and pick $\bv_{t}$ to be some unit vectors orthogonal to $\bw_1,\bw_2,\ldots,\bw_t$, as well as all previously constructed $\bv$'s (this is always possible since the dimension is sufficiently large). By \lemref{lem:info}, as long as all vectors thus constructed are orthogonal to $\bw_t$, the information $\{F(\bw_t),\nabla F(\bw_t),\nabla^2 F(\bw_t)\}$ provided to the algorithm does not depend on $\bv_{t+1},\ldots,\bv_{T}$, and only depends on $\bv_1,\ldots,\bv_t$ which were already determined. Therefore, the next query point $\bw_{t+1}$ is fixed. 
	\item At the end of the process, we pick $\bv_T$ to be some unit vector orthogonal to all previously chosen $\bv$'s as well as $\bw_1,\ldots,\bw_T$. 
\end{itemize}
Based on this construction, the following lemma is self-evident:
\begin{lemma}\label{lem:subopt1}
	It holds that $\inner{\bw_{T},\bv_{T}}=0$.
\end{lemma}

Based on this lemma, we now turn to argue that $\bw_{T}$ must be a sub-optimal point. We first establish the following result:
\begin{lemma}\label{lem:subopt2}
	Letting $\bw^\star=\arg\min_{\bw}F(\bw)$, it holds that
	\[
	\left\|\bw^\star-\sum_{i=1}^{T}q^i \bv_i\right\|~\leq~\sqrt{\frac{T\mu r^2}{16\lambda}}
	\]
	where $q=\frac{\sqrt{\kappa}-1}{\sqrt{\kappa}+1}$.
\end{lemma}
\begin{proof}
	Let $F_r$ denote $F$, where we make the dependence on the parameter $r$ explicit. 
	We first argue that
	\begin{equation}\label{eq:diff}
	\sup_{\bw\in \reals^d} |F_r(\bw)-F_0(\bw)|~\leq \frac{T\mu r^2}{32}.
	\end{equation}
	This is because
	\begin{align*}
	|F_r(\bw)-F_0(\bw)| ~\leq~ \frac{\lambda(\kappa-1)}{8}&\Bigg(\sum_{i=1}^{T-1}\left|\phi_r(\inner{\bv_i-\bv_{i+1},\bw})-\phi_0(\inner{\bv_i-\bv_{i+1},\bw})\right|\\
	&~~~~~+\left|\phi_r(\inner{\bv_T,\bw})-\phi_0(\inner{\bv_T,\bw})\right|\Bigg),
	\end{align*}
	and since $\sup_{z\in\reals}|\phi_r(z)-\phi_0(z)|=\sup_{z\in \reals}|\phi_r(z)-z^2|\leq 2r^2$, the above is at most $\frac{\lambda(\kappa-1)}{4}Tr^2\leq \frac{\lambda\kappa}{4}Tr^2$. Recalling that $\kappa=\mu/8\lambda$, \eqref{eq:diff} follows.
	
	Let $\bw_r=\arg\min F_r(\bw)$. By $\lambda$-strong convexity of $F_0$ and $F_r$,
	\[
	F_0(\bw_r)-F_0(\bw_0) ~\geq~ \frac{\lambda}{2}\norm{\bw_r-\bw_0}^2~~~,~~~F_r(\bw_0)-F_r(\bw_r) ~\geq~ \frac{\lambda}{2}\norm{\bw_0-\bw_r}^2.
	\] 
	Summing the two inequalities and using \eqref{eq:diff},
	\[
	\lambda\norm{\bw_r-\bw_0}^2 ~\leq~ F_0(\bw_r)-F_r(\bw_r)+F_r(\bw_0)-F_0(\bw_0)~\leq~ \frac{T\mu r^2}{16}~,
	\]
	and therefore
	\begin{equation}\label{eq:diff2}
	\norm{\bw_r-\bw_0}^2 ~\leq~ \frac{T\mu r^2}{16\lambda}.
	\end{equation}
	By definition, $\bw_r=\bw^\star$ from the statement of our lemma, so it only remains to prove that $\bw_0=\arg\min F_0(\bw)$ equals $\sum_{i=1}^{T}q^i\bv_i$. To see this, note that $F_0(\bw)$ can be equivalently written as $\tilde{F}(V\bw)$, where $V$ is some orthogonal $d\times d$ matrix with its first $T$ rows equal to $\bv_1,\ldots,\bv_{T}$, and
	\[
	\tilde{F}(\bx) =  \frac{\lambda(\kappa-1)}{8}\left(x_1^2+\sum_{i=1}^{T-1}(x_i-x_{i+1})^2+(a_{\kappa}-1)x_{T}^2-w_1\right)+\frac{\lambda}{2}\norm{\bx}^2.
	\]	
	By an immediate corollary of \lemref{lem:tridig}, $\tilde{F}(\cdot)$ is minimized at $(q,q^2,\ldots,q^{T},0,\ldots,0)$, where $q=\frac{\sqrt{\kappa}-1}{\sqrt{\kappa}+1}$, and therefore $F(\bw)=\tilde{F}(V\bw)$ is minimized at $V^\top(q,q^2,\ldots,q^{T},0,\ldots,0)$, which equals $\sum_{i=1}^{T}q^i \bv_i$ as required.
\end{proof}
Note that this lemma also allows us to bound the norm of $\bw^\star=\arg\min F(\bw)$, since it implies that
\[
\norm{\bw^\star}~\leq~ \left\|\sum_{i=1}^{T}q^i\bv_i\right\|+\sqrt{\frac{T\mu r^2}{16\lambda}},
\]
and since $(a+b)^2\leq 2a^2+2b^2$ and $q<1$, we have
\begin{align*}
\norm{\bw^\star}^2 ~&\leq~ 2\left\|\sum_{i=1}^{T}q^i\bv_i\right\|^2+\frac{T\mu r^2}{8\lambda}~=~ 2\sum_{i=1}^{T}q^{2i}+\frac{T\mu r^2}{8\lambda}\\
&\leq~2\sum_{i=1}^{\infty}q^{2i}+\frac{T\mu r^2}{8\lambda} ~=~ \frac{2q^2}{1-q^2}+\frac{T\mu r^2}{8\lambda}\\
&\leq~ \frac{2}{1-q}+\frac{T\mu r^2}{8\lambda} = \sqrt{\kappa}+1+\frac{T\mu r^2}{8\lambda},
\end{align*}
which is at most $\sqrt{\kappa}+2\leq 3\sqrt{\kappa}$, since we assume that $c$ is sufficiently small so that $\frac{T\mu r^2}{8\lambda}\leq 1$, and that $\kappa=\mu/8\lambda\geq 1$.

The proof of the theorem follows by combining \lemref{lem:subopt1} and \lemref{lem:subopt2}. Specifically, \lemref{lem:subopt1} (which states that $\inner{\bw_T,\bv_{T}}=0$) and the fact that $\bv_1,\ldots,\bv_{T}$ are orthonormal tells us that
\begin{align*}
\left\|\bw_T-\sum_{i=1}^{T}q^i \bv_i\right\|^2 &~=~ \left\|\left(\bw_T-\sum_{i=1}^{T-1}q^i \bv_i\right)-q^{T}\bv_{T}\right\|^2~=~
\left\|\bw_T-\sum_{i=1}^{T-1}q^i \bv_i\right\|^2+\norm{q^T\bv_T}^2\\
&~\geq~ \norm{q^T\bv_T}^2 ~=~ q^{2T},
\end{align*}
and hence 
\[
\left\|\bw_T-\sum_{i=1}^{T}q^i \bv_i\right\| ~\geq q^T.
\]
On the other hand, \lemref{lem:subopt2} states that 
\[
\left\|\bw^\star-\sum_{i=1}^{T}q^i \bv_i\right\|~\leq~\sqrt{\frac{T\mu r^2}{16\lambda}}.
\]
Combining the last two displayed equations by the triangle inequality, we get that
\[
\left\|\bw_T-\bw^\star\right\|~\geq~ q^{T}-\sqrt{\frac{T\mu r^2}{16\lambda}}.
\]
By the assumption that $c$ is sufficiently small so that $\sqrt{\frac{T\mu r^2}{16\lambda}}\leq \frac{1}{2}q^T$, the left hand side is at least $\frac{1}{2}q^{T}$. Squaring both sides, we get
\[
\norm{\bw_T-\bw^\star}^2 ~\geq~ \frac{1}{4}q^{2T},
\]
so by strong convexity of $F$,
\[
F(\bw_T)-F(\bw^\star) ~\geq~ \frac{\lambda}{2}\norm{\bw_T-\bw^\star}^2 ~\geq~ \frac{\lambda}{8}q^{2T}.
\]
Plugging in the value of $q$, we get
\[
F(\bw_T)-F(\bw^\star) ~\geq~  \frac{\lambda}{8}\left(\frac{\sqrt{\kappa}-1}{\sqrt{\kappa}+1}\right)^{2T}.
\]
On the other hand, we showed earlier that $\norm{\bw^\star}^2 \leq 3\sqrt{\kappa}$, so by smoothness, $F(\mathbf{0})-F(\bw^\star) \leq \frac{\mu}{2}\norm{\bw^\star}^2\leq\frac{3\mu}{2}\sqrt{\kappa}$. Therefore,
\[
\frac{F(\bw_T)-F(\bw^\star)}{F(\mathbf{0})-F(\bw^\star)} ~\geq~ \frac{\lambda}{12\mu\sqrt{\kappa}}\left(\frac{\sqrt{\kappa}-1}{\sqrt{\kappa}+1}\right)^{2T}
\]

To make the right-hand side less than $\epsilon$, $T$ must be such that
\[
\left(\frac{\sqrt{\kappa}-1}{\sqrt{\kappa}+1}\right)^{2T} \leq \frac{12\mu\sqrt{\kappa}\epsilon}{\lambda},
\]
which is equivalent to
\[
2T\cdot\log\left(\frac{\sqrt{\kappa}+1}{\sqrt{\kappa}-1}\right) \geq \log\left(\frac{\lambda}{12\mu\sqrt{\kappa}\epsilon}\right).
\]
Since $\log\left(\frac{\sqrt{\kappa}+1}{\sqrt{\kappa}-1}\right) = \log\left(1+\frac{2}{\sqrt{\kappa}-1}\right)\leq \frac{2}{\sqrt{\kappa}-1}$, it follows that $T$ must be such that
\[
\frac{4T}{\sqrt{\kappa}-1}~\geq~ \log\left(\frac{\lambda}{12\mu\sqrt{\kappa}\epsilon}\right).
\] 
Plugging in $\kappa=\mu/8\lambda$ and simplifying a bit, we get that
\[
T ~\geq~ \frac{1}{4}\left(\sqrt{\frac{\mu}{8\lambda}}-1\right)\cdot \log\left(\frac{\sqrt{8}(\lambda/\mu)^{3/2}}{12\epsilon}\right),
\]
from which the result follows.

\subsection{Proof of \thmref{thm:main}}

We will define a randomized choice of quadratic functions $f_1,\ldots,f_n$, and prove a lower bound on the expected optimization error of any algorithm (where the expectation is over both the algorithm and the randomized functions). This implies that for any algorithm, the same lower bound (in expectation over the algorithm only) holds for some deterministic choice of $f_1,\ldots,f_n$. 

There will actually be two separate constructions, one leading to a lower bound of $\Omega(n)$, and one leading to a lower bound of  $\Omega\left(\sqrt{\frac{n\mu}{\lambda }}\cdot\log\left(\frac{(\lambda/\mu)^{3/2}\sqrt{n}}{\epsilon}\right)\right)$. Choosing the construction which leads to the larger lower bound, the theorem follows.

\subsubsection{\texorpdfstring{An $\Omega(n)$ Lower Bound}{An Omega(n) Lower Bound}}

Starting with the $\Omega(n)$ lower bound, let $\delta_i$, where $i\in \{1,\ldots,n\}$, be chosen uniformly at random from $\{-1,+1\}$, and define 
\[
f_{i}(\bw)=- \delta_i w_1+\frac{\lambda}{2}\norm{\bw}^2.
\]
Clearly, these are $\lambda$-smooth (and hence $\mu$-smooth) functions, as well as $\lambda$-strongly convex. Also, the optimum of $F(\bw)=\frac{\mu}{n}\sum_{i=1}^{n}f_i(\bw)$ equals $\bw^\star=\left(\frac{1}{n\lambda}\sum_{i=1}^{n}\delta_i \right)\be_1$, where $\be_1$ is the first unit vector. As a result, $\norm{\bw^\star}^2 = \frac{1}{\lambda^2}\left(\frac{1}{n}\sum_{i=1}^{n}\delta_i\right)^2$, so by $\lambda$-smoothness of $F$
\[
F(\mathbf{0})-F(\bw^\star) ~\leq~ \frac{\lambda}{2}\norm{\bw^\star}^2~=~\frac{1}{2\lambda}\left(\frac{1}{n}\sum_{i=1}^{n}\delta_i\right)^2.
\]
Since $\delta_i$ are i.i.d., we have by Hoeffding's bound that with probability at least $3/4$, $\left|\frac{1}{n}\sum_{i=1}^{n}\delta_i\right|$ is at most $ \sqrt{2\log(8/3)/n}\leq \sqrt{2/n}$. Plugging into the equation above, we get that with probability at least $3/4$,
\begin{equation}\label{eq:main00}
F(\mathbf{0})-F(\bw^\star) ~\leq~\frac{1}{\lambda n}.
\end{equation}

Turning to lower bound $F(\bw_T)-F(\bw^\star)$, we have by strong convexity that
\begin{align*}
F(\bw_T)-F(\bw^\star) &~\geq~ \frac{\lambda}{2}\norm{\bw_T-\bw^\star}^2~\geq~ \frac{\lambda}{2}(w_{T,1}-w^\star_1)^2\\
&~=~
\frac{1}{2\lambda}\left(\lambda w_{T,1}-\frac{1}{n}\sum_{i=1}^{n}{\delta_i}\right)^2.
\end{align*}
Now, if at most $\lfloor n/2 \rfloor$ indices $\{1,\ldots,n\}$ were queried by the algorithm, then the $\bw_T$ returned by the algorithm must be independent of at least $\lceil n/2 \rceil$ random variables $\delta_{j_1},\ldots,\delta_{j_{\lceil n/2\rceil}}$ (for some distinct indices $j_1,j_2,\ldots$ depending on the algorithm's behavior, but independent of the values of $\delta_{j_1},\ldots,\delta_{j_{\lceil n/2 \rceil}}$). Therefore, conditioned on $j_1,\ldots,j_{\lceil n/2\rceil}$ and the values of $\delta_{j_1},\ldots,\delta_{j_{\lceil n/2\rceil}}$, the expression above can be written as
\[
\frac{1}{2\lambda}\left(\eta-\frac{1}{n}\sum_{i\notin\{j_1,\ldots,j_{\lceil n/2\rceil}\}}\delta_i\right)^2~,
\]
where $\eta$ is a fixed quantity independent of the values of $\delta_i$ for $i\notin\{j_1,\ldots,j_{\lceil n/2\rceil}\}$. By a standard anti-concentration argument, with probability at least $3/4$, this expression will be at least $\frac{1}{2\lambda}\left(\frac{c'}{\sqrt{n}}\right)^2=\frac{c'^2}{2\lambda n}$ for some universal positive $c'>0$. Since this is true for any $j_1,\ldots,j_{\lceil n/2\rceil}$ and $\delta_{j_1},\ldots,\delta_{j_{\lceil n/2\rceil}}$, we get that with probability at least $3/4$ over $\delta_1,\ldots,\delta_n$,
\[
F(\bw_T)-F(\bw^\star)~\geq~ \frac{c'^2}{2\lambda n}.
\]
Combining this with \eqref{eq:main00} using a union bound, we have that with probability at least $1/2$,
\[
\frac{F(\bw_T)-F(\bw^\star)}{F(\mathbf{0})-F(\bw^\star)}~\geq~ \frac{c'^2\lambda n}{2\lambda n} ~=~ \frac{c'^2}{2}.
\]
As a result, since the ratio above is always a non-negative quantity,
\[
\E\left[\frac{F(\bw_T)-F(\bw^\star)}{F(\mathbf{0})-F(\bw^\star)}\right]~\geq~ \frac{c'^2}{4}.
\]
Using the assumption stated in the theorem (taking $c=c'^2/4$), we have that the right hand side cannot be smaller than $\epsilon$, unless more than $\lfloor n/2 \rfloor = \Omega(n)$ oracle calls are made.

\subsubsection{\texorpdfstring{An $\Omega\left(\sqrt{\frac{n\mu}{\lambda }}\cdot\log\left(\frac{(\lambda/\mu)^{3/2}\sqrt{n}}{\epsilon}\right)\right)$ Lower Bound}{A linear rate lower bound}}

We now turn to prove the $\Omega\left(\sqrt{\frac{n\mu}{\lambda }}\cdot\log\left(\frac{\lambda}{\epsilon}\right)\right)$ lower bound, using a different function construction: Let $j_1,\ldots,j_{d-1}$ be chosen uniformly and independently at random from $\{1,\ldots,n\}$, and define
\begin{align} \label{def:f_i_proof_thm2}
f_i(\bw) = \frac{\mu-\lambda}{8}\left(\sum_{l=1}^{d-1}\mathbf{1}_{j_l=i}(w_l-w_{l+1})^2+\frac{1}{n}\bigg(w_1^2+(a_{\kappa}-1)w_d^2-w_1\bigg)\right)+\frac{\lambda}{2}\norm{\bw}^2.
\end{align}
where $\mathbf{1}_{A}$ is the indicator of the event $i$. Note that these are all $\lambda$-strongly convex functions, as all terms in their definition are convex in $\bw$, and there is an additional  $\frac{\lambda}{2}\norm{\bw}^2$ term. Moreover, they are also $\mu$-smooth: To see this, note that $\nabla^2 f_i(\bw) \preceq \frac{ (\mu-\lambda)}{4}A+\lambda I\preceq \mu I$, where $A\preceq 4I$ is as defined in the proof of \lemref{lem:tridig}. 

The average function $F(\bw)=\frac{1}{n}\sum_{i=1}^{n}f_i(\bw)$ equals
\begin{align} \label{def:F_proof_thm2}
F(\bw) = \frac{\mu - \lambda }{8n}\left(w_1^2+\sum_{i=1}^{d-1}(w_i-w_{i+1})^2+(a_{\kappa}-1)w_d^2-w_1\right)+\frac{\lambda}{2}\norm{\bw}^2,
\end{align}

Therefore, by \lemref{lem:tridig}, the smoothness parameter of $F$ is $(\mu-\lambda)/n+\lambda\le\mu$, the global minimum $\bw^\star$ of $F$ equals $(q,q^2,\ldots,q^d)$, where $q=\frac{\sqrt{\kappa}-1}{\sqrt{\kappa}+1}$ and $$\kappa=\frac{\frac{\mu-\lambda}{n}+\lambda}{\lambda}=\frac{\frac{\mu}{\lambda}-1}{n}+1.$$ Note that since $q<1$ and $\kappa\geq 1$, the squared norm of $\bw^\star$ is at most
\begin{align}\label{eq:main123}
\sum_{i=1}^{d}q^{2i} ~\leq \sum_{i=1}^{\infty}q^{2i} = \frac{q^2}{1-q^2} \leq \frac{1}{1-q} = \frac{\sqrt{\kappa}+1}{2}~\leq~\sqrt{\kappa},
\end{align}
hence by smoothness,
\begin{equation}\label{eq:main0}
F(\mathbf{0})-F(\bw^\star) ~\leq~ \frac{\mu}{2}\norm{\bw^\star}^2 ~\leq~ \frac{\mu}{2}\sqrt{\kappa}.
\end{equation}

With these preliminaries out of the way, we now turn to compute a lower bound on the expected optimization error. The proof is based on arguing that $\bw_T$ can only have a first few coordinates being non-zero. To see how this gives a lower bound, let $l_T\in \{1,\ldots,d\}$ be the largest index of a non-zero coordinate of $\bw_T$ (or $0$ if $\bw_T=\mathbf{0}$). By definition of $\bw^\star$, we have
\[
\norm{\bw_T-\bw^\star}^2 \geq \sum_{i=l_T+1}^{d} q^{2i} \geq g(l_T),
\] 
where
\begin{align}\label{def:g_proof_thm2}
g(z) = \begin{cases} q^{2(z+1)} & z<d \\ 0 & z\geq d\end{cases}~.
\end{align}
By strong convexity of $F$, this implies that
\[
F(\bw_T)-F(\bw^\star) ~\geq~ \frac{\lambda}{2}\norm{\bw_T-\bw^\star}^2 ~\geq~
\frac{\lambda}{2}g(l_T).
\]
Finally, taking expectation over the randomness of $j_1,\ldots,j_{d-1}$ above (and over the internal randomness of the algorithm, if any), applying \lemref{lem:gq_to_ell}, and choosing the dimension $d=\lceil 2\E[l_T]\rceil$ (which we will later show to equal the value specified in the theorem), we have
\[
\E\left[F(\bw_T)-F(\bw^\star)\right] ~\geq~ \frac{\lambda}{4} q^{4\E[l_T]+4} ~=~ \frac{\lambda}{4}\left(\frac{\sqrt{\kappa}-1}{\sqrt{\kappa}+1}\right)^{2\E[l_T]+2}.
\]
Combined with \eqref{eq:main0}, this gives

\begin{equation}\label{eq:ell_suboptbound}
\E\left[\frac{F(\bw_T)-F(\bw^\star)}{F(\mathbf{0})-F(\bw^\star)}\right]~\geq~ \frac{\lambda}{2\mu\sqrt{\kappa}}\left(\frac{\sqrt{\kappa}-1}{\sqrt{\kappa}+1}\right)^{2\E[l_T]+2}.
\end{equation}
Thus, it remains to upper bound $\E[l_T]$. 

To get a bound, we rely on the following key lemma (where $\be_i$ is the $i$-th unit vector, and recall that $\Wcal_t$ defines the set of allowed query points $\bw_t$, and $j_1,\ldots,j_d$ are the random indices used in constructing $f_1,\ldots,f_n$):
\begin{lemma}\label{lem:lt}
	For all $t$, it holds that $\Wcal_t\subseteq \text{span}\{\be_d,\be_1,\be_2,\be_3,\ldots,\be_{\ell_t}\}$ for all $t$, where $\ell_t$ is defined recursively as follows: $\ell_1=1$, and $\ell_{t+1}$ equals the largest number in $\{1,\ldots,d-1\}$ such that $\{j_{\ell_t},j_{\ell_t+1},\ldots,j_{\ell_{t+1}-1}\}\subseteq\{i_t,i_{t-1},\ldots,i_{\max\{1,t-\lfloor n/2 \rfloor+1\}}\}$ (and $\ell_{t+1}=\ell_t$ if no such number exists).
\end{lemma}
As will be seen later, $\ell_T$ (which is a random variable as a function of the random indices $j_1,\ldots,j_d$) upper-bounds the number of non-zero coordinates of $\bw_T$, and therefore we can upper bound $\E[l_T]$ by $\E[\ell_T]$. 
\begin{proof}
	The proof is by induction over $t$. Since $\Wcal_1=\{\mathbf{0}\}\subseteq \text{span}(\be_d)$, the result trivially holds for $t=1$. Now, suppose that $\Wcal_t\subseteq \text{span}\{\be_d,\be_1,\ldots,\be_{\ell_t}\}$ for some $t$ and $\ell_t$. Note that in particular, this means that $\bw_t$ is non-zero only in its first $\ell_t$ coordinates. By definition of $f_i$ for any $i$,
	\begin{align*}
	\nabla f_{i}(\bw) &= 
	\frac{\lambda n (\kappa-1)}{8}\left(2\sum_{l=1}^{d-1}\mathbf{1}_{j_l=i}(w_l-w_{l+1})(\be_l-\be_{l+1})+\frac{1}{n}\left(2w_1\be_1+2(a_{\kappa}-1)w_d\be_d-\be_1\right)\right)+\lambda \bw\\
	\nabla^2 f_{i}(\bw) &= 
	\frac{\lambda n (\kappa-1)}{8}\left(\sum_{l=1}^{d-1}\mathbf{1}_{j_l=i}(2E_{l,l}-E_{l+1,l}-E_{l,l+1})+\frac{1}{n}\left(2E_{1,1}+2(a_{\kappa}-1)E_{d,d}\right)\right)+\lambda I,	
	\end{align*}
	where $E_{r,s}$ is the $d\times d$ which is all zeros, except for an entry of $1$ in location $(r,s)$. It is easily seen that these expressions imply the following:
	\begin{itemize}
		\item If $j_{\ell_t}\neq i_t$, then $\nabla f_{i_t}(\bw_t)\in \text{span}\{\be_d,\be_1,\ldots,\be_{\ell_t}\}$, otherwise $\nabla f_{i_t}(\bw_t)\in \text{span}\{\be_d,\be_1,\ldots,\be_{\ell_t+1}\}$.
		\item For any $\bw$ and $l\in \{1,\ldots,d-1\}$, if $j_l\neq i$, then $\nabla^2 f_{i}(\bw)$ is block-diagonal, with a block in the first $l\times l$ entries. In other words, any entry $(r,s)$ in the matrix, where $r\leq l$ and $s>l$ (or $r>l$ and $s\leq l$) is zero.
		\item As a result, if $j_l\notin \{i_t,i_{t-1},\ldots,i_{\max\{1,t-\lfloor n/2 \rfloor+1\}}\}$, then $ \sum_{\tau=\max\{1,t-\lfloor n/2 \rfloor+1\}}^{t}\alpha_\tau \nabla^2 f_{i_{\tau}}(\bw_{\tau})$, for arbitrary scalars $\tau$, is block-diagonal with a block in the first $l\times l$ entries. The same clearly holds for any matrix with the same block-diagonal structure. 
	\end{itemize}
	Together, these observations imply that the operations specified in Assumption \ref{assump:linalg} can lead to vectors outside $\text{span}\{\be_d,\be_1,\ldots,\be_{\ell_t}\}$, only if $j_{\ell_t}\in \{i_t,i_{t-1},\ldots,i_{\max\{1,t-\lfloor n/2 \rfloor+1\}}\}$. Moreover, these vectors must belong to $\text{span}\{\be_d,\be_1,\ldots,\be_{\ell_{t+1}}\}$, where $\ell_{t+1}$ is as specified in the lemma: By definition, $j_{\ell_{t+1}}$ is not in $\{i_t,i_{t-1},\ldots,i_{\max\{1,t-\lfloor n/2 \rfloor+1\}}\}$, and therefore all relevant Hessians have a block in the first $\ell_{t+1}\times \ell_{t+1}$ entries, hence it is impossible to create a vector with non-zero coordinates (using the operations of Assumption \ref{assump:linalg}) beyond the first $\ell_{t+1}$. 
\end{proof}

Since $\bw_T\subseteq \Wcal_T$, the lemma above implies that $\E[l_T]$ from \eqref{eq:ell_suboptbound} (where $l_T$ is the largest index of a non-zero coordinate of $\bw_T$) can be upper-bounded by $\E[\ell_T]$, where the expectation is over the random draw of the indices $j_1,\ldots,j_{d-1}$. This can be bounded using the following lemma:
\begin{lemma}
	It holds that $\E[\ell_T]\leq 1+\frac{2(T-1)}{n}$. 
\end{lemma}
\begin{proof}
	By definition of $\ell_t$ and linearity of expectation, we have 
	\begin{equation}\label{eq:lTT}
	\E[\ell_T] = \E\left[\sum_{t=1}^{T-1}\left(\ell_{t+1}-\ell_{t}\right)\right]+\ell_1 = \sum_{t=1}^{T-1}\E[\ell_{t+1}-\ell_{t}]+1.
	\end{equation}
	Let us consider any particular term in the sum above. Since $\ell_{t+1}-\ell_t$ is a non-negative integer, we have
	\[
	E[\ell_{t+1}-\ell_{t}]~=~\Pr\left(\ell_{t+1}>\ell_t\right)\cdot \E\left[\ell_{t+1}-\ell_t ~\middle|~ \ell_{t+1}>\ell_t \right].
	\]
	By definition of $\ell_t$, the event $\ell_{t+1}>\ell_t$ can occur only if $j_{\ell_{t}}\notin \{i_{t-1},i_{t-2},\ldots,i_{\max\{1,t-\lfloor n/2 \rfloor\}}\}$, yet $j_{\ell_{t}}\in \{i_{t},i_{t-1},\ldots,i_{\max\{1,t-\lfloor n/2 \rfloor+1\}}\}$. This is equivalent to $j_{\ell_{t}}=i_{t}$ (that is, in iteration $t$ we happened to choose the index $j_{\ell_t}$ of the unique individual function, which contains the block linking coordinate $\ell_{t}$ and $\ell_{t}+1$, hence allowing us to ``advance'' and have more non-zero coordinates). But since the algorithm is oblivious, $i_{t}$ is fixed whereas $j_{\ell_{t}}$ is chosen uniformly at random, hence the probability of this event is $1/n$. Therefore, $\Pr\left(\ell_{t+1}>\ell_t\right)\leq 1/n$. Turning to the conditional expectation of $\ell_{t+1}-\ell_t$ above, it equals the expected number of indices $j_{\ell_{t}},j_{\ell_{t}+1},\ldots$ belonging to $\{i_{t},i_{t-1},\ldots,i_{\max\{1,t-\lfloor n/2 \rfloor+1\}}\}$, conditioned on $j_{\ell_t}$ belonging to that set. But since the $i$ indices are fixed and the $j$ indices are chosen uniformly at random, this equals one plus the expected number of times where a randomly drawn $j\in \{1,\ldots,n\}$ belongs to  $\{i_{t},i_{t-1},\ldots,i_{t-\lfloor n/2 \rfloor+1}\}$. Since this set contains at most $\lfloor n/2 \rfloor$ distinct elements in $\{1,\ldots,n\}$, this is equivalent to (one plus) the expectation of a geometric random variable, where the success probability is at most $1/2$. By a standard derivation, this is at most $1+\frac{1/2}{1-1/2} = 2$. Plugging into the displayed equation above, we get that 
	\[
	\E[\ell_{t+1}-\ell_t]~\leq~\frac{1}{n}\cdot 2~=~ \frac{2}{n},
	\]
	and therefore the bound in \eqref{eq:lTT} is at most $\frac{2(T-1)}{n}+1$ as required.
\end{proof}
Plugging this bound into \eqref{eq:ell_suboptbound}, we get
\[
\E\left[\frac{F(\bw_T)-F(\bw^\star)}{F(\mathbf{0})-F(\bw^\star)}\right]~\geq~ \frac{\lambda}{2\mu\sqrt{\kappa}}\left(\frac{\sqrt{\kappa}-1}{\sqrt{\kappa}+1}\right)^{\frac{4(T-1)}{n}+4}~.
\]
To make the right-hand side less than $\epsilon$, $T$ must be such that
\[
\left(\frac{\sqrt{\kappa}-1}{\sqrt{\kappa}+1}\right)^{\frac{4(T-1)}{n}+4} \leq \frac{2\mu\sqrt{\kappa}\epsilon}{\lambda},
\]
which is equivalent to
\[
\left(\frac{4(T-1)}{n}+4\right)\log\left(\frac{\sqrt{\kappa}+1}{\sqrt{\kappa}-1}\right) \geq \log\left(\frac{\lambda}{2\mu\sqrt{\kappa}\epsilon}\right).
\]
Since $\log\left(\frac{\sqrt{\kappa}+1}{\sqrt{\kappa}-1}\right) = \log\left(1+\frac{2}{\sqrt{\kappa}-1}\right)\leq \frac{2}{\sqrt{\kappa-1}}$ (see, e.g., Lemma 12 in \cite{arjevani2016dimension}), it follows that $T$ must be such that
\[
\left(\frac{4(T-1)}{n}+4\right)\frac{2}{\sqrt{\kappa-1}}~\geq~ \log\left(\frac{\lambda}{2\mu\sqrt{\kappa}\epsilon}\right).
\] 
Plugging in $\kappa=\frac{\frac{\mu}{\lambda}-1}{n}+1$, we get that
\[
T ~\geq~ 1+\frac{n}{4}\left(\frac{\sqrt{\frac{\frac{\mu}{\lambda}-1}{n}}}{2}\cdot\log\left(\frac{\lambda}{2\mu\epsilon\sqrt{\frac{\frac{\mu}{\lambda}-1}{n}+1}}\right)-4\right).
\]
Using asymptotic notation the right-hand side equals
\[
 \Omega\left(\sqrt{n(\mu/\lambda-1)}\log\left(\frac{(\lambda/\mu)^{3/2}\sqrt{n}}{\epsilon}\right)\right).
\]
as required. The bound on the dimension $d$ follows from the fact that we chose it to be $\Ocal(\E[l_T])=\Ocal(1+T/n)$, and to make the lower bound valid it is enough to pick some $T=\Ocal\left(\sqrt{\frac{n\mu}{\lambda }}\cdot\log\left(\frac{(\lambda/\mu)^{3/2}\sqrt{n}}{\epsilon}\right)\right)$.

													\subsection{Proof of \thmref{thm:mainrelaxed}}

Recall that the proof of \thmref{thm:main} essentially shows that for any (possibly stochastic) index-oblivious optimization algorithm there exists some `bad' assignment of the $d-1$ blocks $j_1,\dots,j_{d-1}$ whose corresponding $f_i:\reals^d\to\reals$  (see \eqref{def:f_i_proof_thm2}) form a functions which is hard-to-optimize. When considering non-oblivious (i.e., adaptive) algorithms this construction fails as soon as the algorithm obtains the Hessians of all the individual functions (potentially, after $n$ second-order oracle queries). Indeed, knowing the Hessians of $f_i$, one can devise an index-schedule which gains no less than one coordinate at every iteration, as opposed to $1/n$ on average for the oblivious case. 
Thus, in order to tackle the non-oblivious case, we form a function over some $D$-dimensional space which `contains' all the $n^{d-1}$ sub-problems at one and the same time (clearly, to carry out our plans we must pick $D$ which grows exponentially fast with $d$, the dimension of the sub-problems). This way, any index-schedule, oblivious or adaptive, must `fit' all the $n^{d-1}$ sub-problems well, and as such, bound to a certain convergence rate which we analyze below. \\

Denote $[n]=\{1,\dots,n\}$, set $D=n^{d-1}d$ and define for any $\bj\in[n]^{d-1}$ the following, 
\begin{align*}
f^\bj_i&: \reals^d\to\reals,&&\bw \mapsto \frac{\mu-\lambda}{8}\left(\sum_{l=1}^{d-1}\mathbf{1}_{j_l=i}(w_l-w_{l+1})^2+\frac{1}{n}\bigg(w_1^2+(a_{\kappa}-1)w_d^2-w_1\bigg)\right)+\frac{\lambda}{2}\norm{\bw}^2,\\
Q^\bj&:\reals^{D}\to\reals^d,&&\bu\mapsto \sum_{l=1}^d \bu^\top\be_{\#\bj d + l}
\end{align*}
where $\#\bj$ enumerates the $n^{d-1}$ tuples $[n]^{d-1}$ from 0 to $n^{d-1}-1$. Note that $f_i^\bj$ are defined exactly as in \eqref{def:f_i_proof_thm2}, only here we make the dependence on $\bj$ explicit. The individual functions are defined as follows:
\begin{align*}
	f_i(\bu)&= \sum_{\bj\in [n]^{d-1}} f^\bj_i( Q^\bj\bu).
\end{align*}
Note that, 
\begin{align*}
\nabla^2 f_i(\bu) 
%&=  \sum_{\bj\in [n]^{d-1}}  \nabla( (Q^\bj)^\top \nabla F( Q^\bj\bu) )\\
%&= \sum_{\bj\in [n]^{d-1}}  (Q^\bj)^\top \nabla( \nabla F( Q^\bj\bu) )\\
&= \sum_{\bj\in [n]^{d-1}} (Q^\bj)^\top \nabla^2 f^\bj_i( Q^\bj\bu) Q^\bj
\end{align*}
Since $\nabla ^2 f_i $ are block-diagonal, we have $\Lambda(\nabla ^2 f_i )=\bigcup_\bj \Lambda(\nabla ^2  f_i^\bj )$, where $\Lambda(\cdot)$ denotes the spectrum of a given matrix. Thus, since $f^\bj_i$ are $\mu$-smooth and $\lambda$-strongly convex (see proof of \thmref{thm:main}), we see that $f_i$ is also $\mu$-smooth and $\lambda$-strongly convex. \\

As for the average function $\Phi(\bu)=\frac{1}{n}\sum_{i=1}^{n}f_i(\bu)$, it is easily verified that for any fixed $\bj\in[n]^{d-1}$,
\begin{align*}
	\frac1n \sum_{i=1}^n f_i^\bj(Q^\bj \bu) = F(Q^\bj \bu), 
\end{align*}
where $F$ is as defined in \eqref{def:F_proof_thm2}. Thus, 
\begin{align*}
\Phi(\bu)=\frac{1}{n} \sum_{i=1}^{n} \sum_{\bj\in [n]^{d-1}} f^\bj_i( Q^\bj\bu)
 =  \sum_{\bj\in [n]^{d-1}}  F( Q^\bj\bu).
\end{align*}
The compute the minimizer of $\Phi$, we compute the first-order derivative:
\begin{align*}
\nabla \Phi(\bu) &= \nabla \left(\sum_{\bj\in [n]^{d-1}}  F( Q^\bj\bu)\right)\\
&=  \sum_{\bj\in [n]^{d-1}} \nabla\bigg( F( Q^\bj\bu)\bigg)\\
&=  \sum_{\bj\in [n]^{d-1}} (Q^\bj)^\top \nabla F( Q^\bj\bu) 
\end{align*}
Thus, by setting $\bu^*=\sum_{\bj} (Q^\bj)^\top \bw^*$, where $\bw^*$ is the minimizer of $F$ as in \lemref{lem:tridig}, we get
\begin{align*}
	\nabla \Phi(\bu^*)&=  \sum_{\bj\in [n]^{d-1}} (Q^\bj)^\top \nabla F\bigg( Q^\bj\sum_{\bj} (Q^\bj)^\top \bu^*\bigg) 
	=  \sum_{\bj\in [n]^{d-1}} (Q^\bj)^\top \nabla F( \bw^*) =0
\end{align*}
Note that, by \eqref{eq:main0}, $\norm{\bu^*}^2=n^{d-1}\norm{\bw^*}^2\le n^{d-1}\sqrt{\kappa}$. Hence, by smoothness,
\begin{equation}\label{def:non_obl_eq1}
\Phi(\mathbf{0})-\Phi(\bu^\star) ~\leq~ \frac{\mu}{2}\norm{\bu^\star}^2 ~\leq~ \frac{\mu}{2}n^{d-1}\sqrt{\kappa}.
\end{equation}

To derive the analytical properties of $\Phi$, we compute the second derivative:
\begin{align*}
\nabla^2 \Phi(\bu) &=  \sum_{\bj\in [n]^{d-1}}  \nabla( (Q^\bj)^\top \nabla F( Q^\bj\bu) )\\
&= \sum_{\bj\in [n]^{d-1}}  (Q^\bj)^\top \nabla( \nabla F( Q^\bj\bu) )\\
&= \sum_{\bj\in [n]^{d-1}} (Q^\bj)^\top \nabla^2 F( Q^\bj\bu) Q^\bj
\end{align*}
Since $\nabla ^2 \Phi $ is a block-diagonal matrix, we have $\Lambda(\nabla ^2 \Phi )=\bigcup_\bj \Lambda(\nabla ^2  F )=\Lambda(\nabla ^2  F )$. Thus, by \lemref{lem:tridig}, it follows that $\Phi$ is $((\mu-\lambda)/n + \lambda)$-smooth and $\lambda$-strongly convex. \\

With these preliminaries out of the way, we now turn to compute a lower bound on the expected optimization error. The proof is based on arguing that $\bu_T$ can only have a first few coordinates being non-zero for each of the $n^{d-1}$ sub-problems. To see how this gives a lower bound, let $l^\bj_T\in \{1,\ldots,d\}$ be the largest index of a non-zero coordinate of $Q^\bj\bu_T$ (or $0$ if $Q^\bj\bu_T=\mathbf{0}$). By definition of $\bu^\star$ and \eqref{eq:main123}, we have
\begin{align*}
	 \norm{\bu_T-\bu^*}^2&= \norm{\sum_{\bj} (Q^\bj)^\top Q^\bj  \bu_T- \sum_{\bj} (Q^\bj)^\top \bw^*}^2\\
	&= \norm{\sum_{\bj} (Q^\bj)^\top( Q^\bj  \bu_T -\bw^*)}^2\\
	&= \sum_{\bj} \norm{ Q^\bj  \bu_T -\bw^*}^2\\
	&\ge \sum_{\bj} g(l^\bj_T)
\end{align*}
where $g$ is defined in \eqref{def:g_proof_thm2}. By the strong convexity of $F$, this implies that
\[
\Phi(\bu_T)-\Phi(\bu^\star) ~\geq~ \frac{\lambda}{2}\norm{\bu_T-\bu^\star}^2 ~\geq~
\frac{\lambda}{2}\sum_{\bj} g(l^\bj_T).
\]

We now proceed along the same lines as in the proof of \thmref{thm:main}. First, to upper bound $l^\bj_T$ (note that, $g$ is monotonically decreasing), we use the following generalized version of \lemref{lem:lt} (whose proof is a straightforward adaptation of the proof of \lemref{lem:lt}):
\begin{lemma}\label{lem:lt_comb}
	Under Assumption \ref{assump:linalg}, for all $t$, it holds that 
\begin{align*}
		\Ucal_t\subseteq \text{span}\left\{\bigcup_{\bj\in[n]^{d-1}}\{\be_{\#\bj d+ d},\be_{\#\bj d+1},\be_{\#\bj d+2},\be_{\#\bj d+3},\ldots,\be_{\#\bj d+\ell^\bj_t}\}\right\}
\end{align*}
	 for all $t$, where $\ell^\bj_t$ is defined recursively as follows: $\ell^\bj_1=1$, and $\ell^\bj_{t+1}$ equals the largest number in $\{1,\ldots,d-1\}$ such that $\{j_{\ell^\bj_t},j_{\ell^\bj_t+1},\ldots,j_{\ell^\bj_{t+1}-1}\}\subseteq\{i_t,i_{t-1},\ldots,i_{\max\{1,t-\lfloor n/2 \rfloor+1\}}\}$ (and $\ell^\bj_{t+1}=\ell^\bj_t$ if no such number exists).
\end{lemma}

As in the proof of \thmref{thm:main}, $\ell_T^\bj$ bound $l_T^\bj$ from above (for any given choice of $i_1,\dots,i_T$), and since $d$ is chosen so that 
\begin{align} \label{ineq:d_thm3}
\frac{1}{n^{d-1}} \sum_{\bj} \ell^\bj_T \le \frac{d}2,
\end{align}
we may take expectation over the internal randomness of the algorithm (if any), and combine it with (\ref{def:non_obl_eq1}) and \lemref{lem:gq_to_ell_comb_thm3} and \lemref{lem:average_comb} below to get
\begin{align*}
\E\left[\frac{\Phi(\bu_T)-\Phi(\bu^\star)}{\Phi(\mathbf{0})-\Phi(\bu^\star)}\right]&\geq
\E\left[\frac{\lambda}{\mu\sqrt{\kappa}n^{d-1}} \sum_{\bj} g(l^\bj_T) \right]
\geq\E\left[\frac{\lambda}{\mu\sqrt{\kappa}n^{d-1}} \sum_{\bj} g(\ell^\bj_T) \right]\\
&\geq \E\left[\frac{\lambda}{2\mu\sqrt{\kappa}} g\left(\frac{1}{n^{d-1}}\sum_{\bj} \ell^\bj_T\right) \right]
\geq \frac{\lambda}{2\mu\sqrt{\kappa}}\left(\frac{\sqrt{\kappa}-1}{\sqrt{\kappa}+1}\right)^{\frac{4(T-1)}{n}+4}.
\end{align*}
Following the same derivation as in the proof of \thmref{thm:main}, we get that $T$ must be of order of 
\[
 \Omega\left(\sqrt{n(\mu/\lambda-1)}\log\left(\frac{(\lambda/\mu)^{3/2}\sqrt{n}}{\epsilon}\right)\right),
\]
as required. The bound on $d$ follows from the fact that we chose it to satisfy Inequality (\ref{ineq:d_thm3}) through the following condition,
\begin{align*}
2\left(1+ \frac{2(T-1)}{n}\right) \le d,
\end{align*} 
and to make the lower bound valid it is enough to pick some 
$T=\Ocal\left(\sqrt{\frac{n\mu}{\lambda 
}}\cdot\log\left(\frac{(\lambda/\mu)^{3/2}\sqrt{n}}{\epsilon}\right)\right)$. 
Thus, we have that $d$ is  $\tilde{\Ocal}(1+\sqrt{\mu/\lambda n})$, implying 
$D=n^{d-1}d=n^{\tilde{\Ocal}\left({1+\sqrt{\mu/\lambda n}}\right)}$.

\begin{lemma}\label{lem:average_comb}
For any fixed sequence $\bi \coloneqq i_1,\dots,i_T\in[n]$ of individual functions chosen during a particular execution of an optimization algorithm which satisfies Assumption 2, it holds that,
\begin{align*}
 \frac{1}{n^{d-1}}\sum_{\bj} \ell_T^\bj \le 1+ \frac{2(T-1)}{n}.
\end{align*}
\end{lemma}

\begin{proof}
By \lemref{lem:lt_comb}, $\ell^\bj_{t+1}$ depends only on $j_p$ for $\ell_t^\bj\le p\le\ell_{t+1}^\bj$. Thus, we may define
\begin{align*}
	A_s = \bigg|\bigg\{(j_1,\dots,j_s)~|~ \ell^{(j_1,\dots,j_s,*)}_{t}=s,~  \ell^{(j_1,\dots,j_s,*)}_{t+1} >s\bigg\}\bigg|,~~~s\in[d],\\B_s = \bigg|\bigg\{(j_1,\dots,j_s)~|~ \ell^{(j_1,\dots,j_s,*)}_{t}=s,~  \ell^{(j_1,\dots,j_s,*)}_{t+1} =s\bigg\}\bigg|,~~~s\in[d].
\end{align*}
Intuitively, $A_s$ and $B_s$ count how many tuples $(j_1,\dots,j_s)$, under a given choice of $i_1,\dots,i_T$, allow at most $s$ non-zero coordinates after $t$ iterations, with one major difference: in $A_s$ we want to allow the algorithm to make a progress after $t+1$ iterations (equivalently, $j_s=i_t$), whereas in $B_s$ we want the algorithm to have the same number of $s$ non-zero coordinates after $t+1$ (equivalently, $j_s\neq i_t$). One can easily verify the following:
\begin{align*}
	&\sum_{s=1}^d (A_s+B_s)n^{d-s-1} = n^{d-1}, \\
	&B_s= (n-1)A_s.
\end{align*}
The first equality may be obtained by splitting the space of all $[n]^{d-1}$ tuples into a group of disjoint sets characterized by the maximal number of non-zero coordinates the algorithm may gain by the $t$ iteration. The second equality is a simple consequence of the way $j_s$ is being constrained by $A_s$ and $B_s$.
This yields,
\begin{align} \label{eq:mainrelax12}
\sum_{s=1}^d A^sn^{-s} = n^{-1}.
\end{align}
Denoting $\Ical\coloneqq\{i_t,i_{t-1},\dots,i_{\max\{t-\lfloor n/2\rfloor+1,1\}}\}$, we get that for any $1\le s\le d-1$ and $1\le k \le d-s$, 
\begin{align*}
	\bigg|\bigg\{\bj~|~ &\ell^\bj_t = s,~ \ell^\bj_{t+1}=s+k  \bigg\}\bigg| \\
	&=\bigg|\bigg\{(j_1,\dots,j_{s-1})~|~ \ell^{(j_1,\dots,j_{s-1},i_t,*)}_{t}=s  \bigg\}\bigg|\cdot
	\bigg|\bigg\{(j_{s+1},\dots,j_{s+k})~|~  j_{s+1},\dots,j_{s+k-1}\in\Ical,~ 
	j_{s+k}\notin\Ical \bigg\}\bigg|\cdot n^{d-s-k-1}\\
	&=A_s |\Ical|^{k-1}(n-|\Ical|) n^{d-s-k-1}
\end{align*}
This allows us to bound from above the average $\ell_{t+1}^\bj - \ell_t^\bj$ over $\bj$ as follows,
\begin{align*}
	\frac1{n^{d-1}}\sum_{\bj}(\ell_{t+1}^\bj-\ell_t^\bj)&=
	\frac1{n^{d-1}}\sum_{s=1}^{d-1} \sum_{k=1}^{d-s}
	\bigg|\bigg\{\bj~|~ \ell^\bj_t = s,~ \ell^\bj_{t+1}=s+k  \bigg\}\bigg|k\\
	&=\frac1{n^{d-1}}\sum_{s=1}^{d-1} \sum_{k=1}^{d-s} A_s |\Ical|^{k-1}(n-|\Ical|) n^{d-s-k-1}k	\\
	&=\sum_{s=1}^{d-1}  A_s n^{-s} \sum_{k=1}^{d-s} |\Ical|^{k-1}(n-|\Ical|) n^{-k}k	\\
	&=\sum_{s=1}^{d-1}  A_s n^{-s} \left(1-\frac{|\Ical|}{n}\right)\sum_{k=1}^{d-s} \left(\frac{|\Ical|}{n}\right)^{k-1} k	\\
	&=\sum_{s=1}^{d-1}  A_s n^{-s} \left(1-\frac{|\Ical|}{n}\right)\sum_{k=1}^{\infty} \left(\frac{|\Ical|}{n}\right)^{k-1} k.
\end{align*}
By standard manipulations of power series we have,
\begin{align*}
	\sum_{k=0}^\infty x^k = \frac1{1-x}  \implies 	\sum_{k=1}^\infty k x^{k-1} = \frac1{(1-x)^2}.
\end{align*}
Combining this with \eqref{eq:mainrelax12} and the fact that $|\Ical|\le n/2$ yields,
\begin{align*}
	\frac1{n^{d-1}}\sum_{\bj}(\ell_{t+1}^\bj-\ell_t^\bj)
	\le\sum_{s=1}^{d-1}  A_s n^{-s} \left(1-\frac{|\Ical|}{n}\right)^{-1}	
	\le2\sum_{s=1}^{d-1}  A_s n^{-s} 
	\le \frac2n,
\end{align*}
which, in turn, gives 
\begin{align*}
	\frac{1}{n^{d-1}}\sum_{\bj} \ell_T^\bj&= \frac{1}{n^{d-1}}\sum_{\bj} \left(\sum_{t=1}^{T-1} (\ell_{t+1}^\bj - \ell_t^\bj)  + \ell_1^\bj \right)\\
	&= \sum_{t=1}^{T-1} \frac{1}{n^{d-1}}\sum_{\bj} (\ell_{t+1}^\bj - \ell_t^\bj) + \frac{1}{n^{d-1}}\sum_{\bj} \ell_1^\bj \\
	&\le \frac{2(T-1)}{n} + 1.
\end{align*}

\end{proof}

%\bibliographystyle{plainnat}
%\bibliography{../mybib}

\begin{lemma}\label{lem:gq_to_ell_comb_thm3} 
	For some $q\in (0,1)$ and positive $d$, define \[g(z) = \begin{cases}q^{2(z+1)}& z< d\\0 & z\geq d\end{cases}~.\] Let $a_1,\dots,a_p$ be a sequence of non-negative reals, such that 
	\begin{align*}
		\frac1p\sum_{i=1} ^p a_i \le \frac{d}2
	\end{align*}
	then 
	\begin{align*}
		\frac1p\sum_{i}^p g(a_i) &\ge \frac12 g\left( \frac1p\sum_{i=1} ^p a_i \right).
	\end{align*}

\end{lemma}
\begin{proof}
	Since $q\in (0,1)$, the function $z\mapsto q^z$ is convex for non-negative $z$. Therefore, by definition of $g$ and Jensen's inequality we have
\begin{align*}
	\frac1p\sum_{i}^p q(a_i) &= \frac{|\{i:a_i<d\}|}p\frac{1}{|\{i:a_i<d\}|}\sum_{\{i:a_i<d\}} g(a_i)\\
	&\ge\frac{|\{i:a_i<d\}|}p g\left(\frac{1}{|\{i:a_i<d\}|}\sum_{\{i:a_i<d\}}a_i\right)\\
\end{align*}
Note that, 
	\begin{align*}
	\frac{d}2 \ge \frac1p\sum_{i=1} ^p a_i = \frac1p\sum_{\{i:a_i<d\}}a_i + \frac1p \sum_{\{i:a_i\ge d\}} a_i\ge
	\frac{d}{p} |\{ i| a_i \ge d\}|\implies
	 	\frac{|\{ i| a_i < d\}|}{p}\ge \frac{1}2
\end{align*}
Therefore, together with the fact that $g$ decreases monotonically and that 
\begin{align*}
	\frac{1}{|\{i:a_i<d\}|}\sum_{\{i:a_i<d\}}a_i\le \frac1p\sum_{i=1} ^p a_i 
\end{align*}
we get 
\begin{align*}
	\frac1p\sum_{i}^p q(a_i) &\ge \frac12 g\left( \frac1p\sum_{i=1} ^p a_i \right)\\
\end{align*}
\end{proof}

\end{document}